\batchmode
\makeatletter
\def\input@path{{D:/Schlumberger/LaTeX/Notes/Bandwidth_Extention/SIAM_Annual_2014/Fractional_Gaussian_Arxiv/}}
\makeatother
\documentclass[oneside,english]{amsart}
\usepackage[T1]{fontenc}
\usepackage[latin9]{inputenc}
\usepackage{amstext}
\usepackage{amsthm}
\usepackage{amssymb}
\usepackage{graphicx}

\makeatletter

\providecommand{\tabularnewline}{\\}

\numberwithin{equation}{section}
\numberwithin{figure}{section}

\pdfoutput=1

\makeatother

\usepackage{babel}
\begin{document}

\title[Fractional derivative of Gaussian and Dawson's integral]{Approximating fractional derivative of the Gaussian function and
Dawson's integral}

\author{Can Evren Yarman}

\address{Schlumberger, High Cross, Madingley Road, Cambridge CB3 0EL, United
Kingdom (cyarman@slb.com)}
\begin{abstract}
A new method for approximating fractional derivatives of the Gaussian
function and Dawson's integral are presented. Unlike previous approaches,
which are dominantly based on some discretization of Riemann-Liouville
integral using polynomial or discrete Fourier basis, we take an alternative
approach which is based on expressing the Riemann-Liouville definition
of the fractional integral for the semi-infinite axis in terms of
a moment problem. As a result, fractional derivatives of the Gaussian
function and Dawson's integral are expressed as a weighted sum of
complex scaled Gaussian and Dawson's integral. Error bounds for the
approximation are provided. Another distinct feature of the proposed
method compared to the previous approaches, it can be extended to
approximate partial derivative with respect to the order of the fractional
derivative which may be used in PDE constraint optimization problems.
\end{abstract}

\maketitle

\section{Introduction}

Fractional derivatives find applications in many branches of physics
and engineering ranging from quantum optics to astrophysics and cosmology,
dynamics of materials to biophysics and medicine, dynamical chaos
to control, signal processing to communications and more. Increase
in the application areas within the last decades \cite{kilbas2006theory}
fueled investigation of numerical methods to compute fractional derivatives.
For recent comprehensive reviews on fractional derivatives and their
applications we refer the reader to \cite{Li20130037,uchaikin2013fractional,kilbas2006theory}
and references within. 

Majority of the numerical methods for computation of fractional derivatives
are based on polynomial expansions, finite difference operators or
discretization of Riemann-Liouville integral over a finite interval
\cite{uchaikin2013fractional,Fukunaga20120152}. For the semi-infinite
case, finite or infinite impulse response filters are designed to
approximate the fractional derivative provided that the functions
can be approximated using discrete convolution or discrete Fourier
transforms\cite{tseng2000computation,tseng2001design,Farid2004,chen2003new,tseng2010design,tseng2013design}.
In this work we took an alternative approach to compute the fractional
derivative of Gaussian function and its Hilbert transform, Dawson's
integral. Our method is based on expressing the Riemann-Liouville
definition of the fractional integral for the semi-infinite axis in
terms of a moment problem. Error bounds for the approximation are
derived. Another advantage of the proposed method to previous approaches,
it can be extended to approximate partial derivatives with respect
to the order of the fractional derivative which may be used in PDE
constraint optimization.

The outline of the paper is as follows. Section \ref{sec:Fractional-integral-and-derivative}
provides background on fractional integral and derivative operators.
In Section \ref{sec:Gaussian-function,-Dawson's-integral}, the method
for approximating fractional derivatives of Gaussian function and
Dawson's integral and corresponding error are derived. Finally in
Section \ref{sec:Partial-a-of-frac_der_a}, approximating partial
derivative with respect to the order of fractional derivative is presented.
Appendix \ref{sec:A-rational-approximation2dawson} provides a rational
approximation to the Dawson's integral which is used in our numerical
results presented in Sections \ref{sec:Gaussian-function,-Dawson's-integral}
and \ref{sec:Partial-a-of-frac_der_a}.

\section{Fractional integral and derivative operators\label{sec:Fractional-integral-and-derivative}}

Fractional derivative of an analytic function $f\left(t\right)$
can be thought as a generalization of Cauchy's integral formula 
\begin{align*}
D^{n}f\left(z\right)=f^{\left(n\right)}\left(z\right) & =\frac{n!}{2\pi\mathrm{i}}\int_{C}\frac{f\left(t\right)}{\left(t-z\right)^{n+1}}dt,
\end{align*}
 from integer $n$ to a real number $\nu$ using generalization of
factorial by Gamma function, $n!=\Gamma\left(n+1\right)$, 
\begin{align*}
f^{\left(\nu\right)}\left(z\right) & =\frac{\Gamma\left(\nu+1\right)}{2\pi\mathrm{i}}\int_{C}\frac{f\left(t\right)}{\left(t-z\right)^{\nu+1}}dt,
\end{align*}
 for some appropriate contour $C$ in the complex plane. Selection
of the contour is important because for non integer $\nu$ the integrand
would contain a branch point rather than a pole. Instead of a closed
contour, which was considered by Sonin and Letnikov, choosing the
contour to be open, Laurent produced the Riemann-Liouville definition
of the fractional integral \cite{miller1993introduction}
\begin{align*}
_{c}D_{x}^{-v}f\left(x\right) & =\frac{1}{\Gamma\left(\nu\right)}\int_{c}^{x}\left(x-t\right)^{\nu-1}f\left(t\right)dt,\quad\text{Re}\left\{ \nu\right\} >0,\,\left(c,x\right)\subset\mathbb{R}.
\end{align*}
Formally, one cannot replace $-\nu$ with $\nu$ to obtain the fractional
derivative operator $_{c}D_{x}^{\nu}$ from the fractional integral
operator because the integral would be divergent. However, writing
$\nu=n-\mu$, where $n=\left\lceil \nu\right\rceil $ is the smallest
integer larger than $\nu$, by analytic continuation, differentiation
for arbitrary order can be defined by$_{c}D_{x}^{\nu}={}_{c}D_{x}^{n}{}_{c}D_{x}^{-\mu}$. 

For $c=0$, the fractional derivative and integral operator satisfies
\begin{align*}
_{0}D_{x}^{\nu}x^{a} & =\frac{\Gamma\left(a+1\right)}{\Gamma\left(a-\nu+1\right)}x^{a-\nu},\quad a\ge0,\nu\in\mathbb{R},
\end{align*}
 and for $c=-\infty$
\begin{align}
_{-\infty}D_{x}^{\nu}\mathrm{e}^{ax} & =a^{\nu}\mathrm{e}^{ax}.\label{eq:fractional_derivative_exponential}
\end{align}
(\ref{eq:fractional_derivative_exponential}) is Leibniz's definition
of fractional derivative (see 4.8.1. in \cite{uchaikin2013fractional}).
When $a$ is imaginary it is also referred to as the Fourier definition
of fractional derivative (see 4.8.3. in \cite{uchaikin2013fractional}).
In computation of the fractional derivative of the Gaussian function
and Dawson's integral, we will consider the case for $c=-\infty$
and, for simplicity of notation, use $d_{x}^{\nu}$ to denote $_{-\infty}D_{x}^{\nu}$.

In the literature $d_{x}^{\nu}$ has been implemented using discrete
convolutions and discrete Fourier transforms through design of finite
or infinite impulse response filters \cite{tseng2000computation,tseng2001design,Farid2004,chen2003new,tseng2010design,tseng2013design}.
We take an alternative approach and get an analytic approximation
to the fractional derivative of the Gaussian and Dawson's integral
by exploiting their Hilbert transform relationship . Because fractional
derivative is a linear operator, if a function can be approximated
as a sum of Gaussian functions (for example as in \cite{mazia2007approximate})
and/or Dawson's integrals, then our method can be used to approximate
its fractional derivative. The proposed method also enables analytic
approximation of derivative with respect to the order of the fractional
derivative of the functions.

\section{Gaussian function, Dawson's integral and approximating their fractional
derivatives\label{sec:Gaussian-function,-Dawson's-integral}}

\subsection{Gaussian function and Dawson's integral}

Dawson's integral $F\left(x\right)$ is defined by (see 7.2(ii) of
\cite{olver2010nist})
\begin{align}
F\left(x\right) & =\frac{1}{2}\int_{0}^{\infty}\mathrm{e}^{-t^{2}/4}\sin\left(xt\right)dt.\label{eq:dawson-definition}
\end{align}
It is related to the Gaussian function $\mathrm{e}^{-t^{2}}$ through
Hilbert transform (see pg 465, (3.5) in \cite{king2009hilbert}):
\begin{align*}
\mathcal{H}\left[\mathrm{e}^{-t^{2}}\right] & =\frac{2}{\sqrt{\pi}}F\left(t\right).
\end{align*}
We define the function $g\left(t\right)=\mathrm{e}^{-t^{2}}+\mathrm{i}\,2\pi^{-1/2}F\left(t\right)$,
whose Taylor series expansion at $t=0$ is given by 
\begin{align}
g\left(t\right) & =\mathrm{e}^{-t^{2}}+\mathrm{i}\frac{2}{\sqrt{\pi}}F\left(t\right)=\sum_{n=0}^{\infty}\frac{\mathrm{i}^{n}}{\Gamma\left(\frac{n+2}{2}\right)}t^{n}.\label{eq:g(t)}
\end{align}
We approximate the fractional derivative $f_{a,\sigma}\left(t\right)=d_{t}^{a}f_{0,\sigma}\left(t\right)$
of $f_{0,\sigma}\left(t\right)=G_{0,\sigma}\left(t\right)+\mathrm{i}\mathcal{H}\left[G_{0,\sigma}\right]\left(t\right)$,
where 
\begin{align*}
G_{0,\sigma}\left(t\right) & =\mathrm{e}^{-\frac{t^{2}}{2\sigma^{2}}}=\sqrt{\frac{2}{\pi}}\int_{0}^{\infty}\mathrm{e}^{-\frac{\omega^{2}}{2}}\cos\left(\frac{\omega}{\sigma}t\right)d\omega\\
\mathcal{H}\left[G_{0,\sigma}\right]\left(t\right) & =\frac{2}{\sqrt{\pi}}F\left(\frac{t}{\sqrt{2}\sigma}\right)=\sqrt{\frac{2}{\pi}}\int_{0}^{\infty}\mathrm{e}^{-\frac{\omega^{2}}{2}}\sin\left(\frac{\omega}{\sigma}t\right)d\omega,
\end{align*}
 as a sum of $g\left(t\right)$:
\begin{align*}
f_{a,\sigma}\left(t\right) & =\sum_{m=1}^{M}\alpha_{m}g\left(\gamma_{m}t\right)+\epsilon\left(t\right),
\end{align*}
 for some $\left(\alpha_{m},\gamma_{m}\right)\in\mathbb{C}^{2}$. 

Consequently, $f_{a,\sigma}\left(t\right)=G_{a,\sigma}\left(t\right)+\mathrm{i}\mathcal{H}\left[G_{a,\sigma}\right]\left(t\right)$,
where $G_{a,\sigma}\left(t\right)=\mathrm{Re}\left\{ f_{a,\sigma}\right\} \left(t\right)=d_{t}^{a}G_{0,\sigma}\left(t\right)$
and $\mathcal{H}\left[G_{a,\sigma}\right]\left(t\right)=\mathrm{Im}\text{\ensuremath{\left\{  f_{a,\sigma}\right\} } }\left(t\right)=d_{t}^{a}\mathcal{H}\left[G_{0,\sigma}\right]\left(t\right)$
are related to the $a^{\text{th}}$ fractional derivative of the Gaussian
function and Dawson's integral, respectively.

\subsection{Approximating fractional derivatives of Gaussian function and Dawson's
integral \label{subsec:Approximating-fractional-derivative}}

Starting with 
\begin{align*}
f_{0,\sigma}\left(t\right) & =\sqrt{\frac{2}{\pi}}\int_{0}^{\infty}\mathrm{e}^{-\frac{\omega^{2}}{2}}\exp\left(\mathrm{i}\frac{\omega}{\sigma}t\right)d\omega,
\end{align*}
Taylor series expansion of $f_{a,\sigma}\left(t\right)$ at $t=0$
is given by {\small{}}{\small \par}

{\small{}
\begin{align}
f_{a,\sigma}\left(t\right) & =\sqrt{\frac{2}{\pi}}\int_{0}^{\infty}\mathrm{e}^{-\frac{\omega^{2}}{2}}\left(\frac{\omega}{\sigma}\right)^{a}\exp\left(\mathrm{i}\left[\frac{\omega}{\sigma}t+a\frac{\pi}{2}\right]\right)d\omega\label{eq:f_a_sigma_fourier}\\
 & =\sqrt{\frac{2}{\pi}}\exp\left(\mathrm{i}a\frac{\pi}{2}\right)\int_{0}^{\infty}\mathrm{e}^{-\frac{\omega^{2}}{2}}\left(\frac{\omega}{\sigma}\right)^{a}\exp\left(\mathrm{i}\frac{\omega}{\sigma}t\right)d\omega\nonumber \\
 & =\sqrt{\frac{2}{\pi}}\exp\left(\mathrm{i}a\frac{\pi}{2}\right)\sum_{n=0}^{\infty}\left[\int_{0}^{\infty}\mathrm{e}^{-\frac{\omega^{2}}{2}}\left(\frac{\omega}{\sigma}\right)^{a}\frac{\left(\mathrm{i}\frac{\omega}{\sigma}\right)^{n}}{n!}d\omega\right]t^{n}\nonumber \\
 & =\sqrt{\frac{2}{\pi}}\exp\left(\mathrm{i}a\frac{\pi}{2}\right)\sum_{n=0}^{\infty}\frac{1}{n!}\frac{\mathrm{i}^{n}}{\sigma^{n+a}}\left[\int_{0}^{\infty}\mathrm{e}^{-\frac{\omega^{2}}{2}}\omega^{a+n}d\omega\right]t^{n}\nonumber \\
 & =\sqrt{\frac{2}{\pi}}\exp\left(\mathrm{i}a\frac{\pi}{2}\right)\sum_{n=0}^{\infty}\frac{1}{n!}\frac{\mathrm{i}^{n}}{\sigma^{n+a}}2^{\left(a+n-1\right)/2}\Gamma\left(\frac{a+n+1}{2}\right)t^{n}.\label{eq:f_a_sigma_taylor}
\end{align}
Equating coefficients of the Taylor series of left and right hand
side of }
\begin{align*}
f_{a,\sigma}\left(t\right) & =\sum_{m=1}^{M}\alpha_{m}g\left(\gamma_{m}t\right)+\epsilon\left(t\right),
\end{align*}
 we obtain the following moment problem for $\left(\alpha_{m},\gamma_{m}\right)$:
\begin{align}
\frac{1}{\sqrt{\pi}}\exp\left(\mathrm{i}a\frac{\pi}{2}\right)\left(\frac{\sqrt{2}}{\sigma}\right)^{a+n}\frac{\Gamma\left(\frac{n+2}{2}\right)\Gamma\left(\frac{a+n+1}{2}\right)}{n!} & =\sum_{m}\alpha_{m}\gamma_{m}^{n}+\varepsilon_{n},\label{eq:frac_derivative-moment_problem}
\end{align}
which can be solved using the method in \cite{YF2014}. Here $\epsilon\left(t\right)=\sum_{n=0}^{\infty}\frac{\epsilon_{n}}{n!}t^{n}$
and $\varepsilon_{n}=\mathrm{i}^{-n}\Gamma\left(\frac{n+2}{2}\right)\frac{\epsilon_{n}}{n!}$.
For $a=2^{-1}$and $\sigma=\sqrt{2}^{-1}$, a solution for the moment
problem (\ref{eq:frac_derivative-moment_problem}) is presented in
Table \ref{tab:(alpha_m,gamma_m)_frac_gaussian_in_gaussian+dawson}. 

\begin{figure}
\begin{tabular}{|c|c|c|}
\hline 
\includegraphics[width=0.3\textwidth]{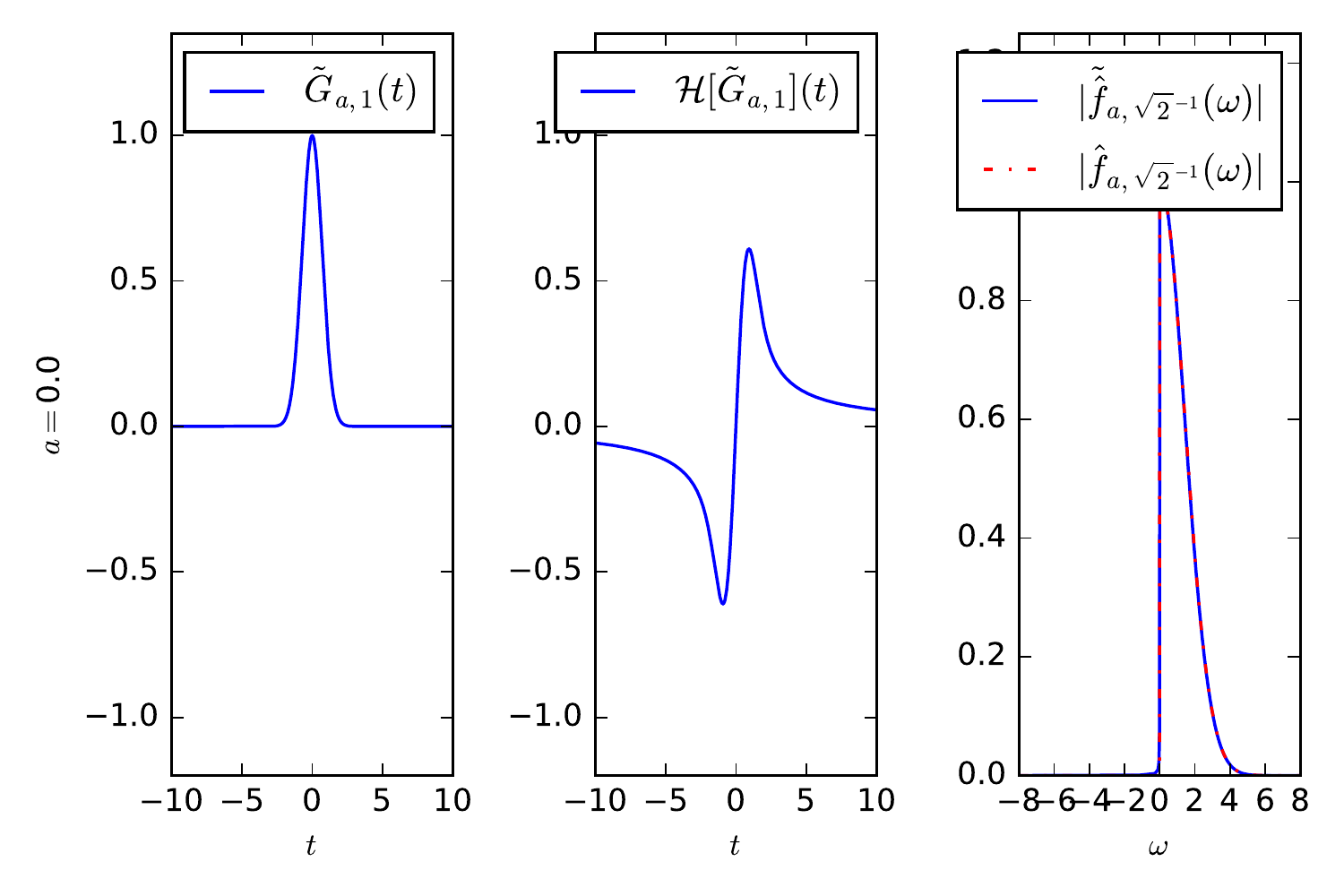} &
\includegraphics[width=0.3\textwidth]{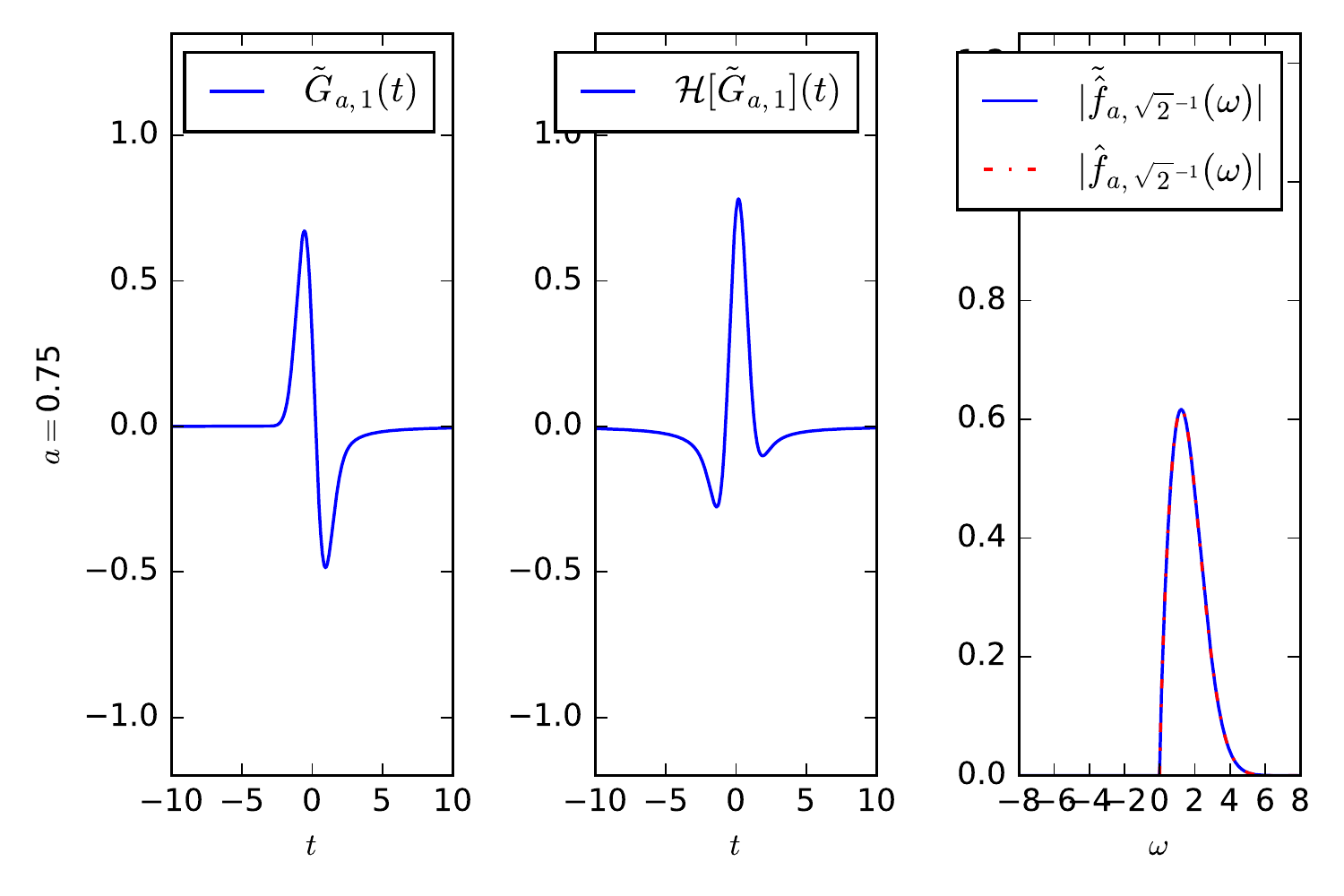} &
\includegraphics[width=0.3\textwidth]{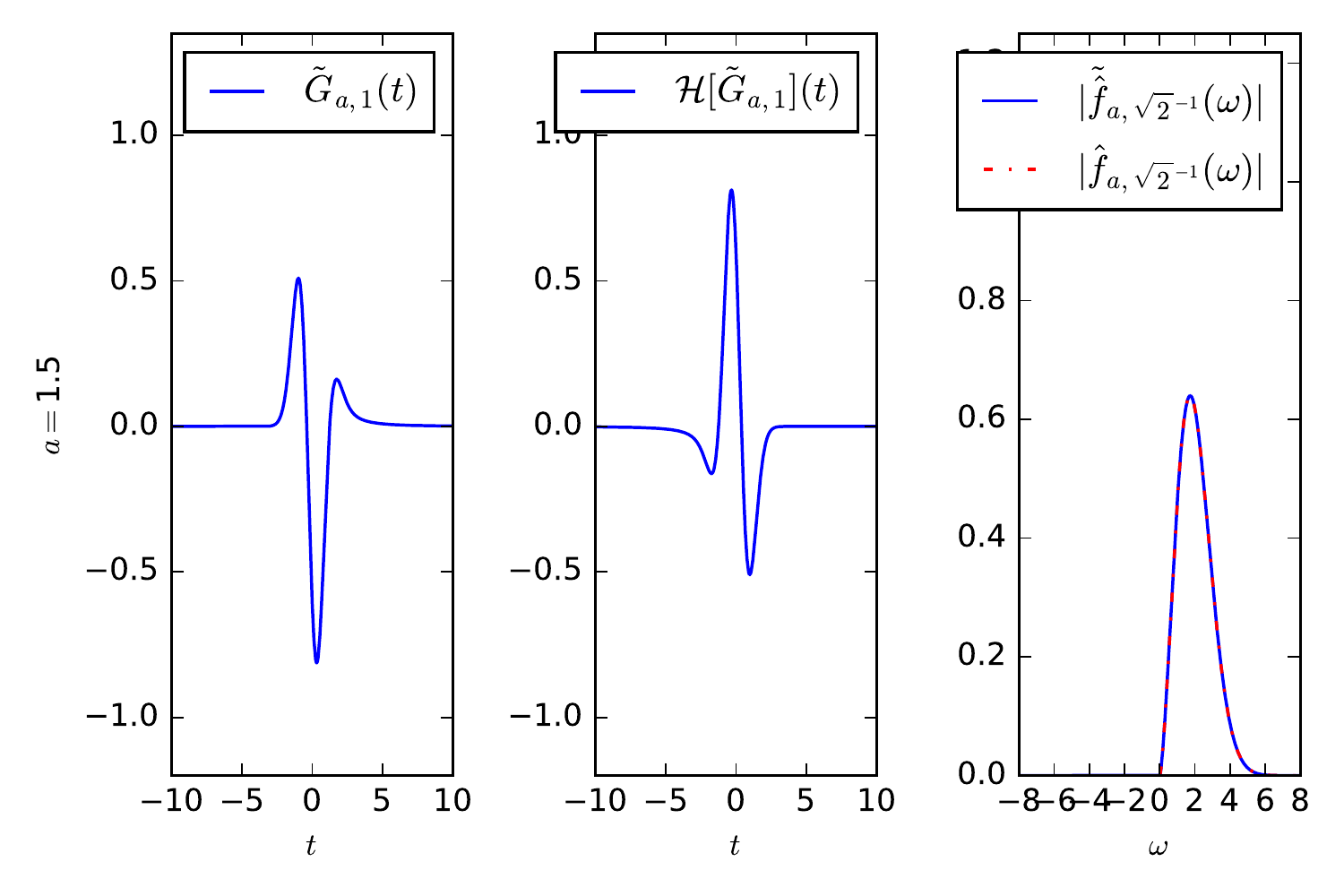}\tabularnewline
\hline 
\includegraphics[width=0.3\textwidth]{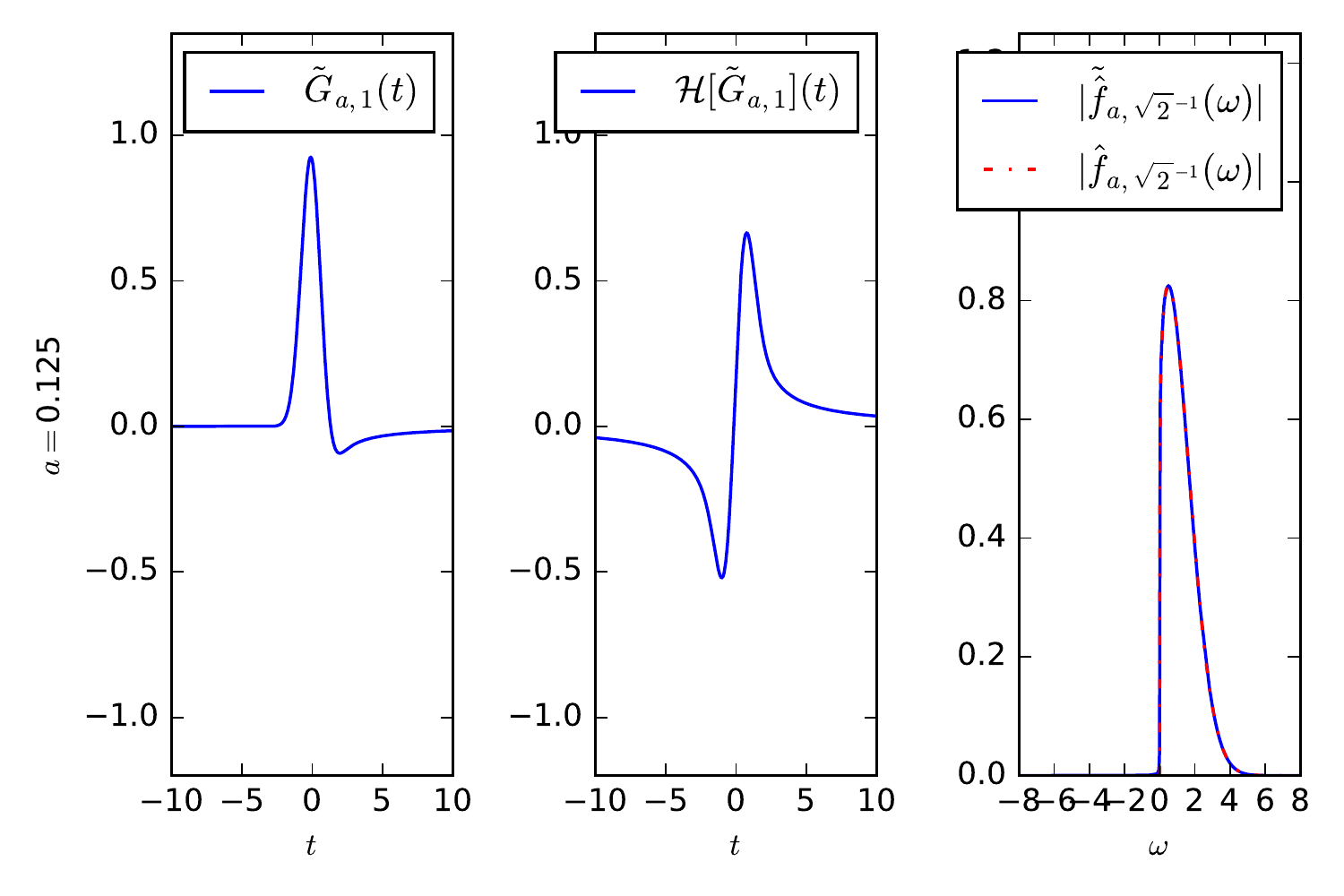} &
\includegraphics[width=0.3\textwidth]{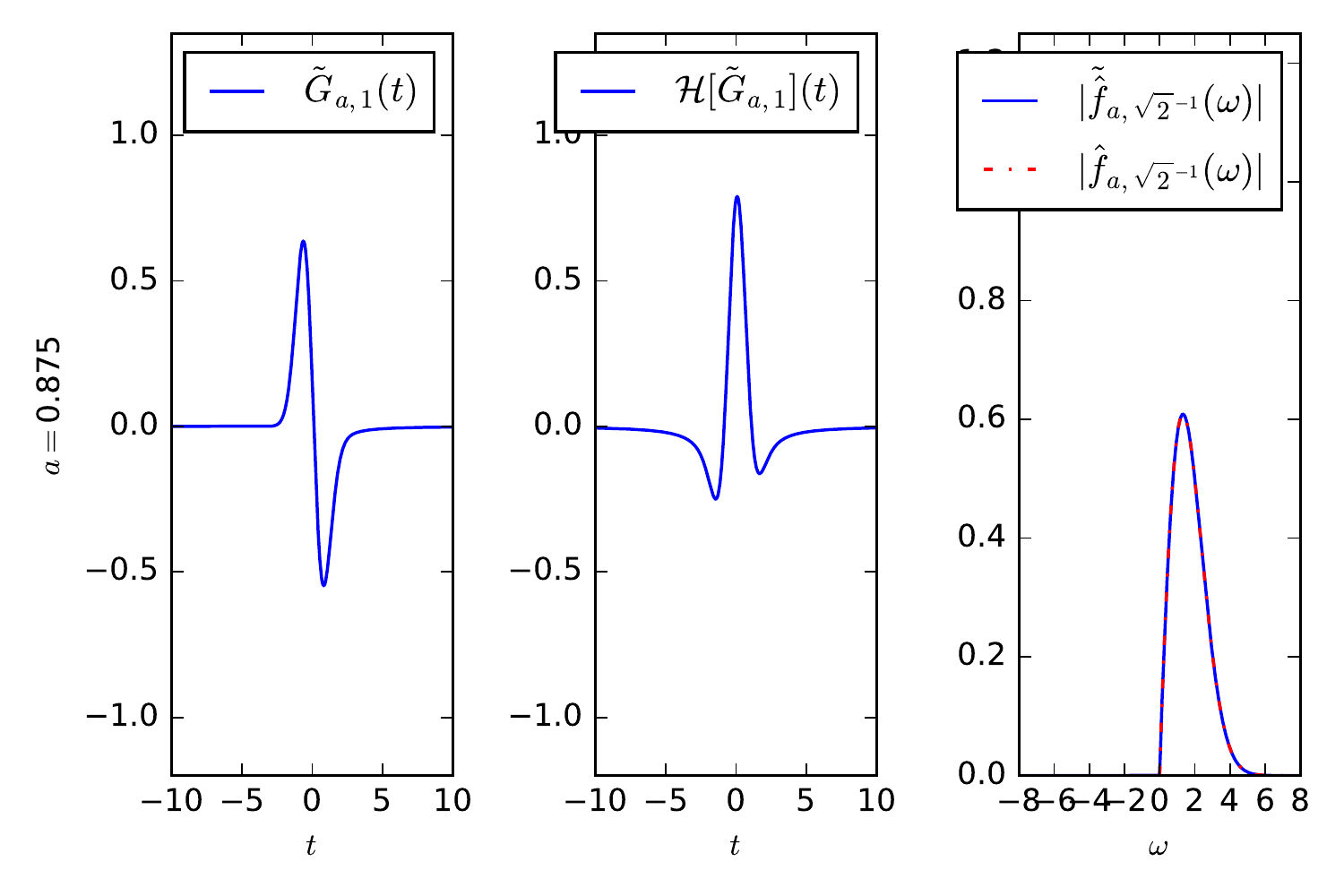} &
\includegraphics[width=0.3\textwidth]{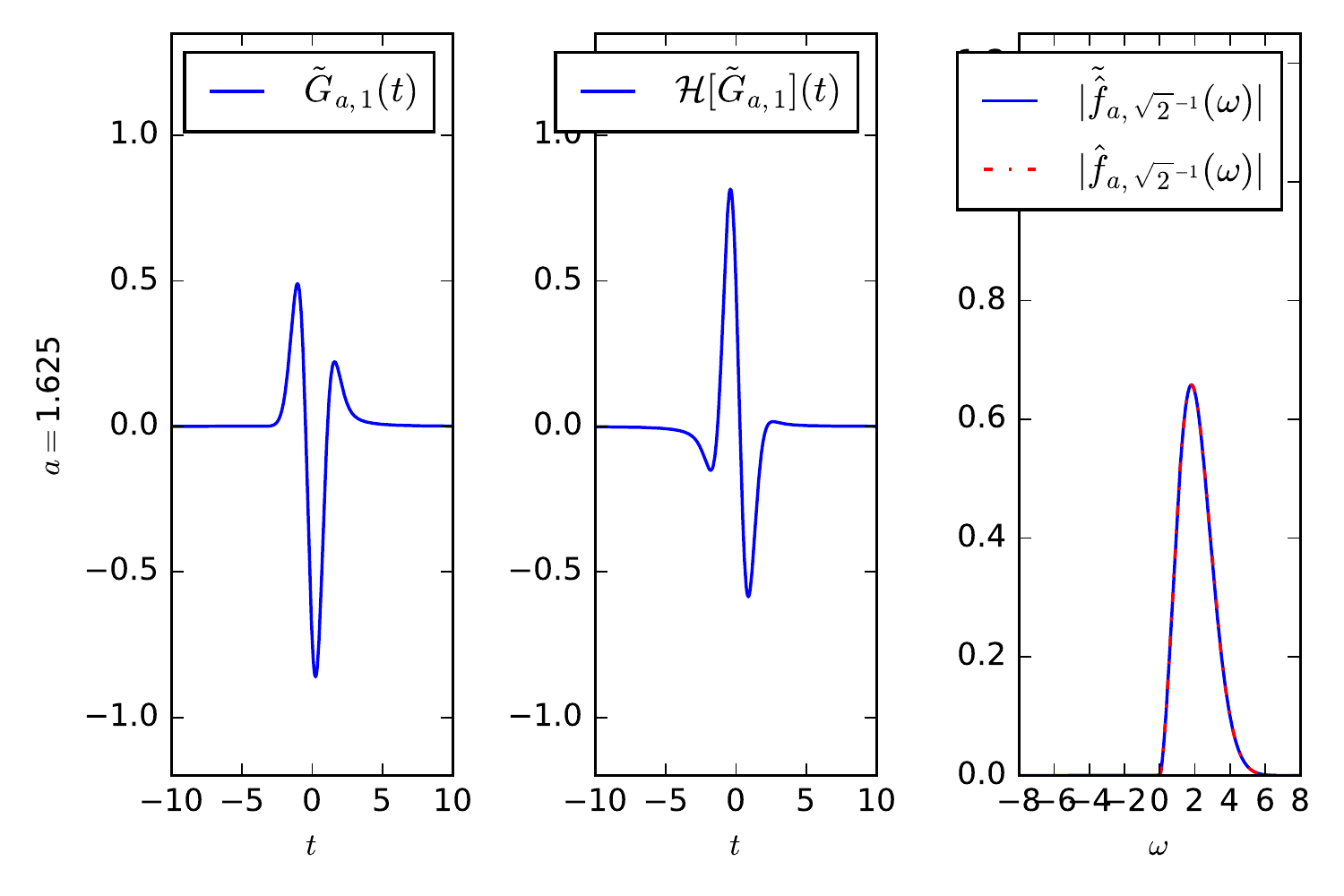}\tabularnewline
\hline 
\includegraphics[width=0.3\textwidth]{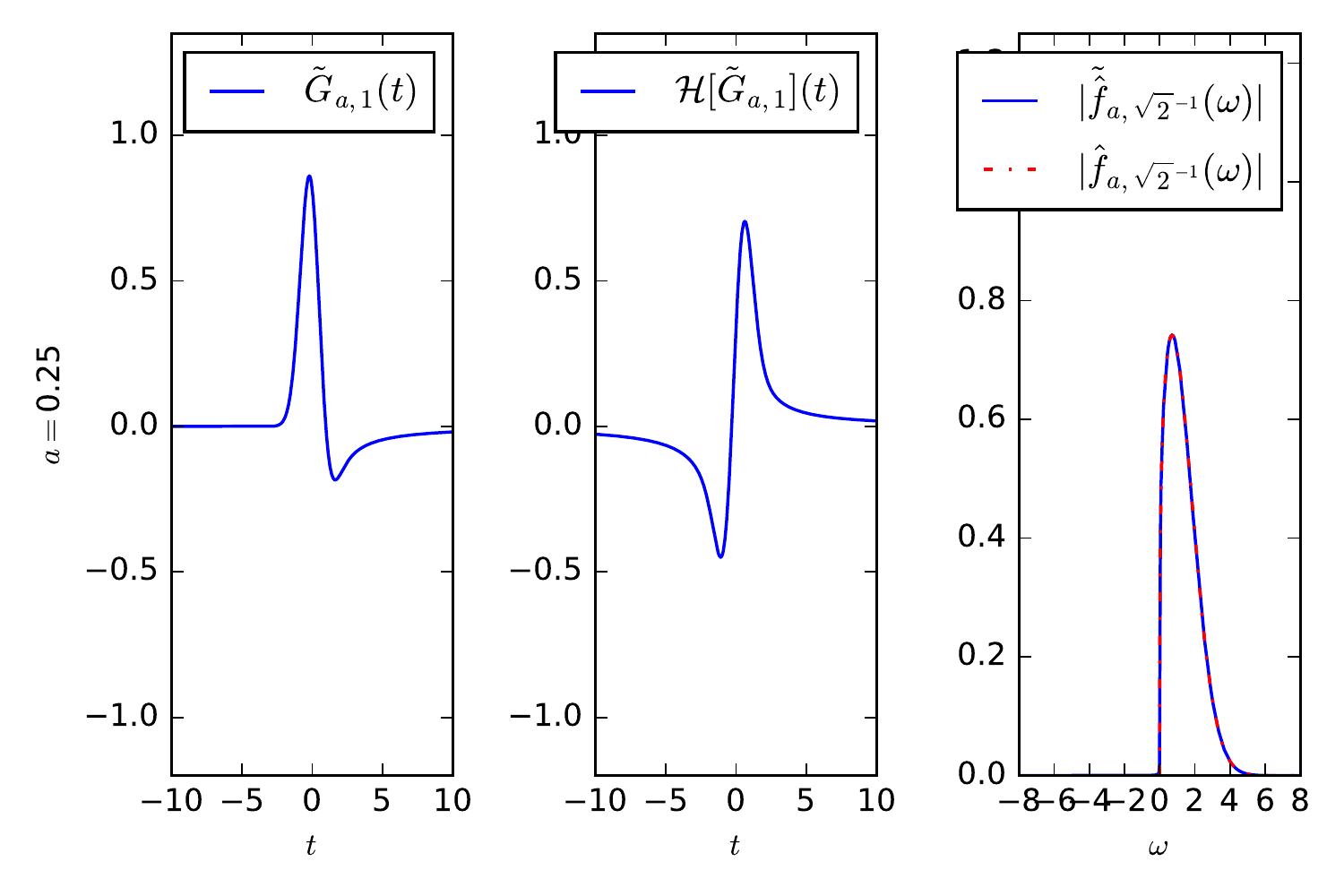} &
\includegraphics[width=0.3\textwidth]{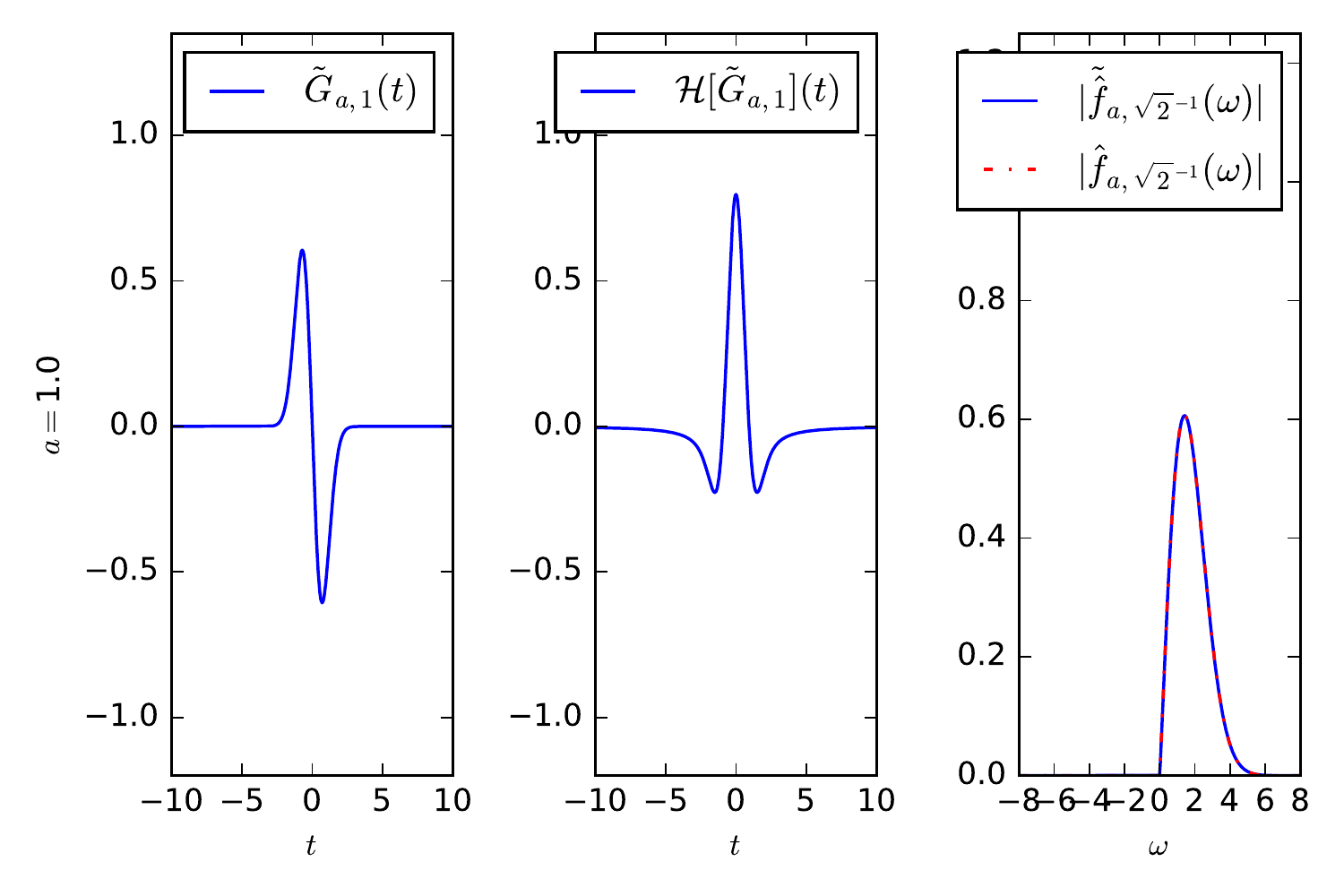} &
\includegraphics[width=0.3\textwidth]{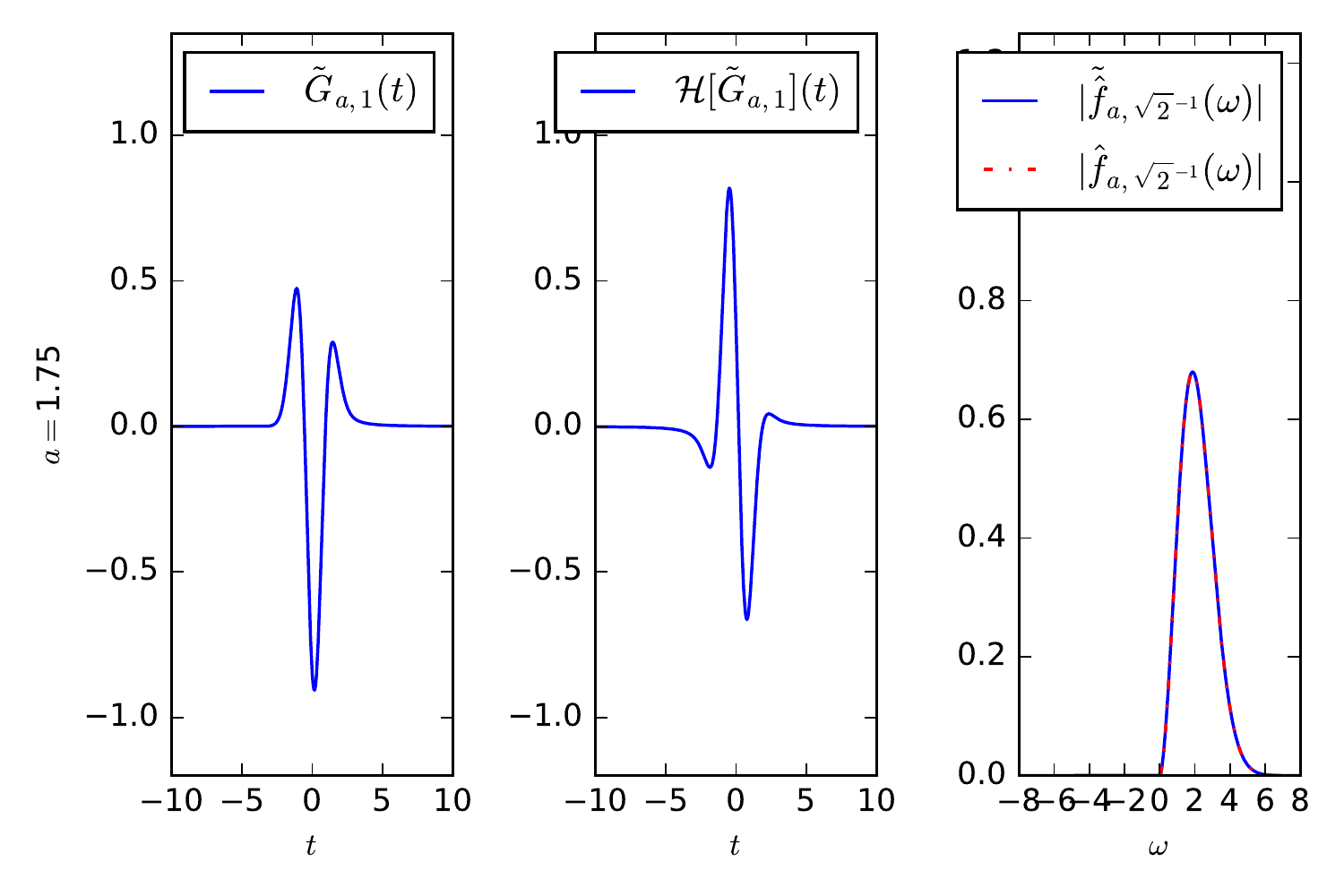}\tabularnewline
\hline 
\includegraphics[width=0.3\textwidth]{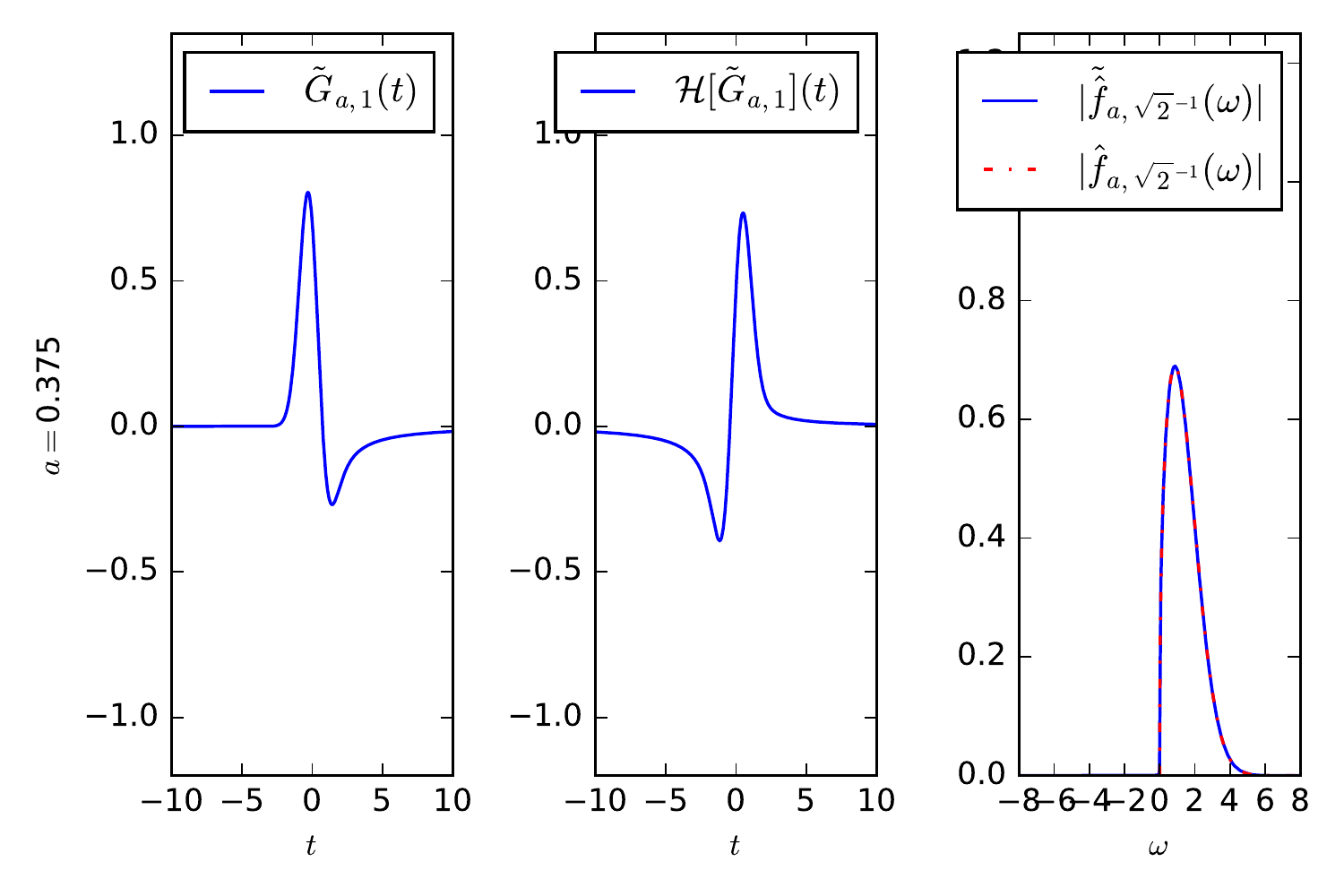} &
\includegraphics[width=0.3\textwidth]{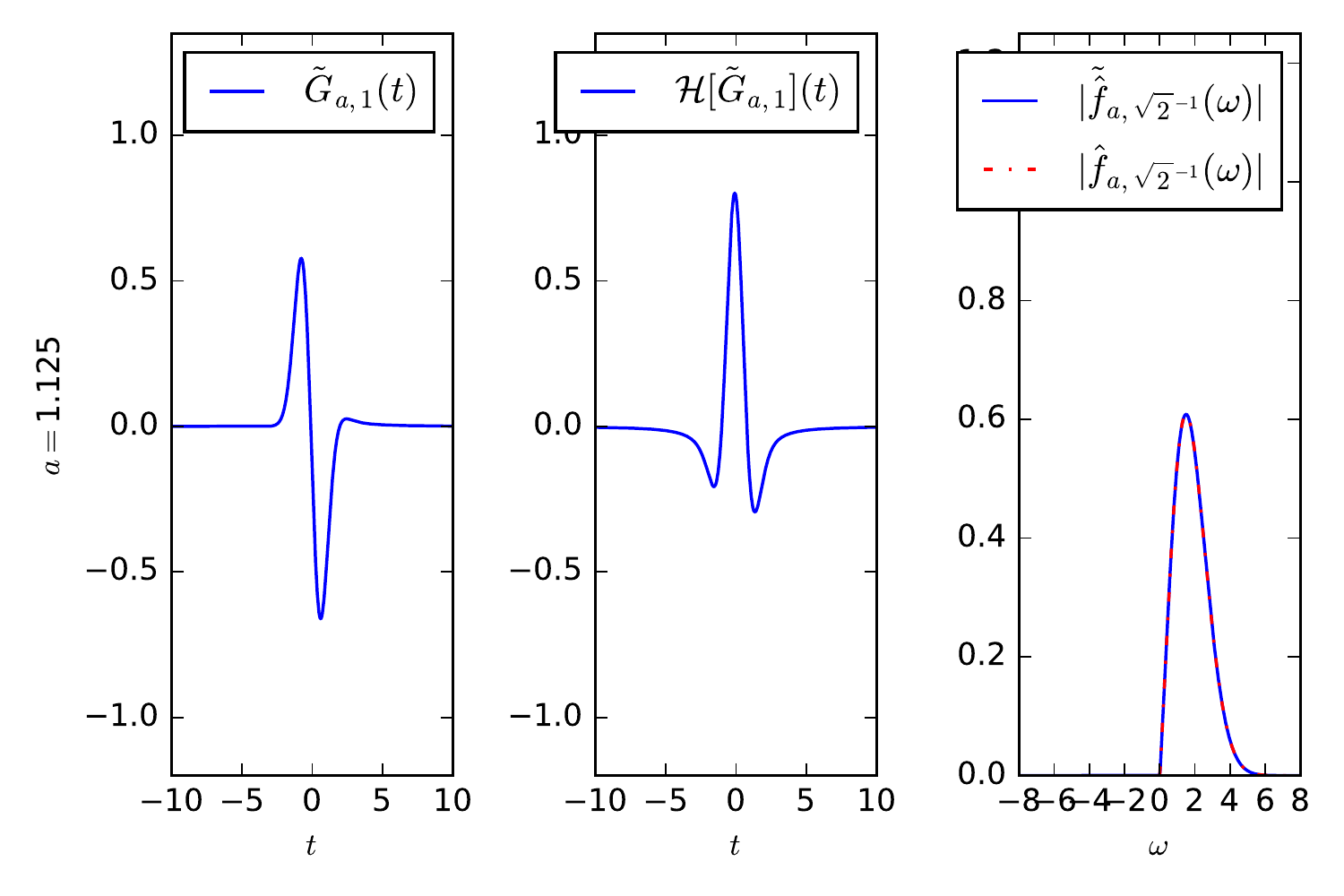} &
\includegraphics[width=0.3\textwidth]{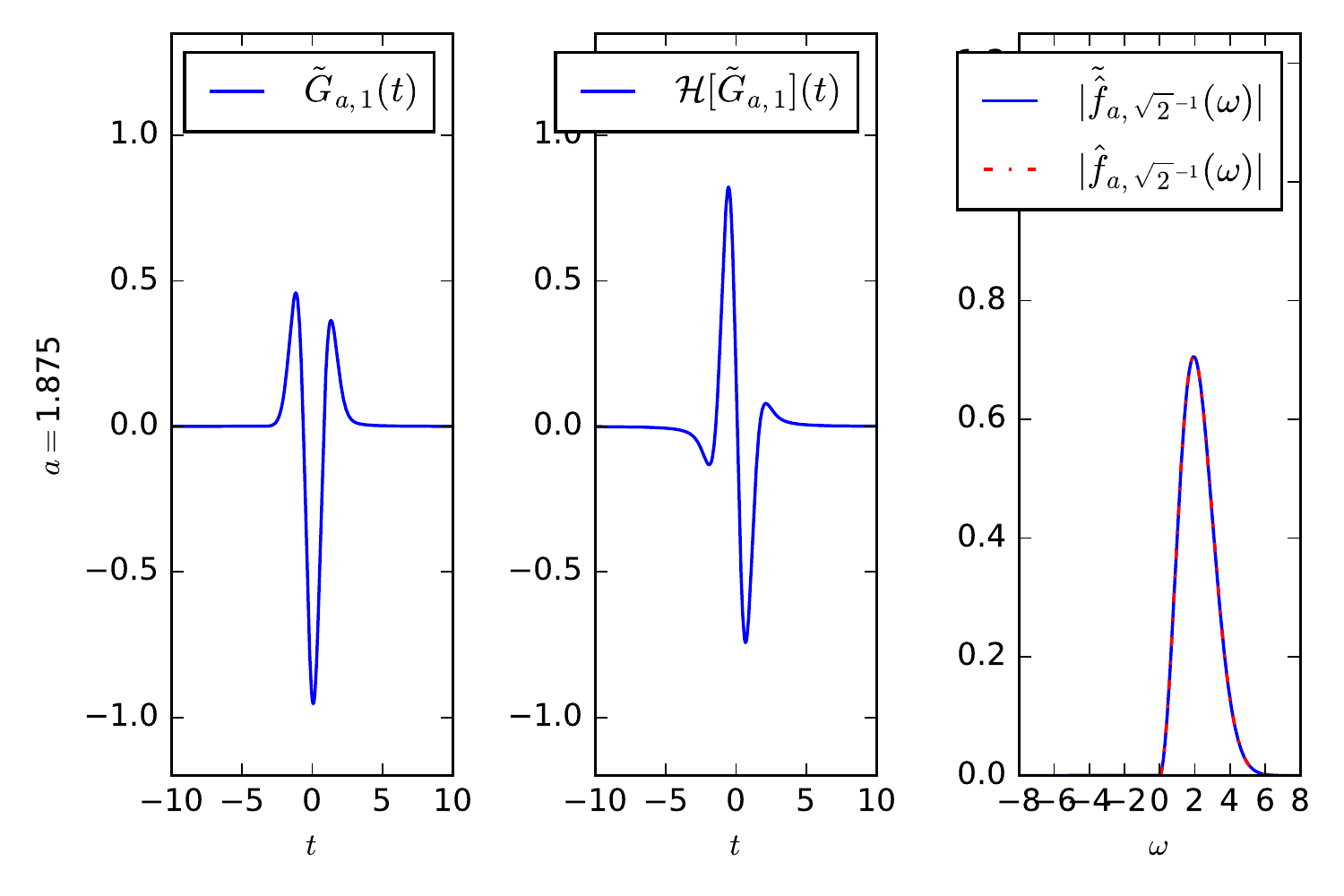}\tabularnewline
\hline 
\includegraphics[width=0.3\textwidth]{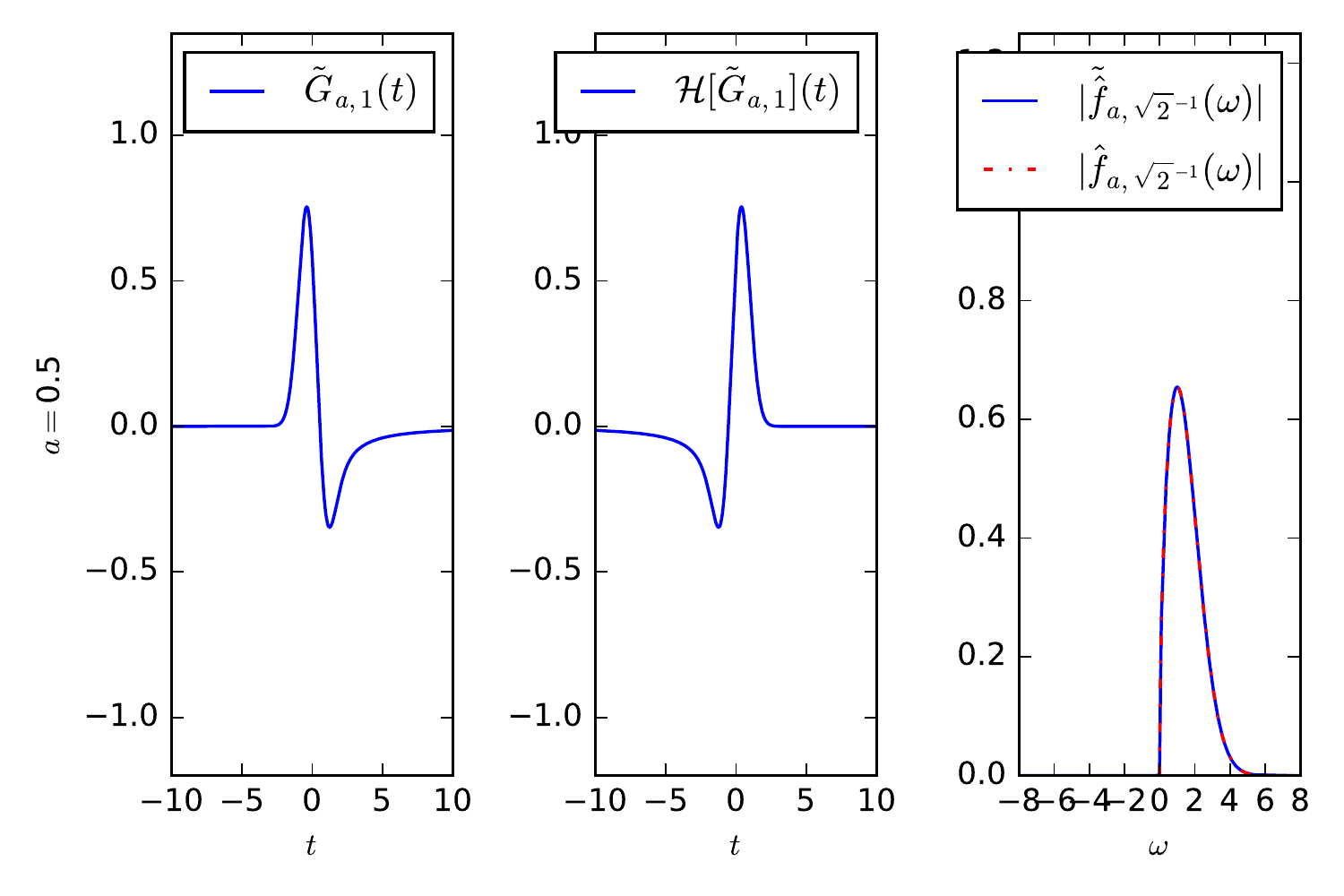} &
\includegraphics[width=0.3\textwidth]{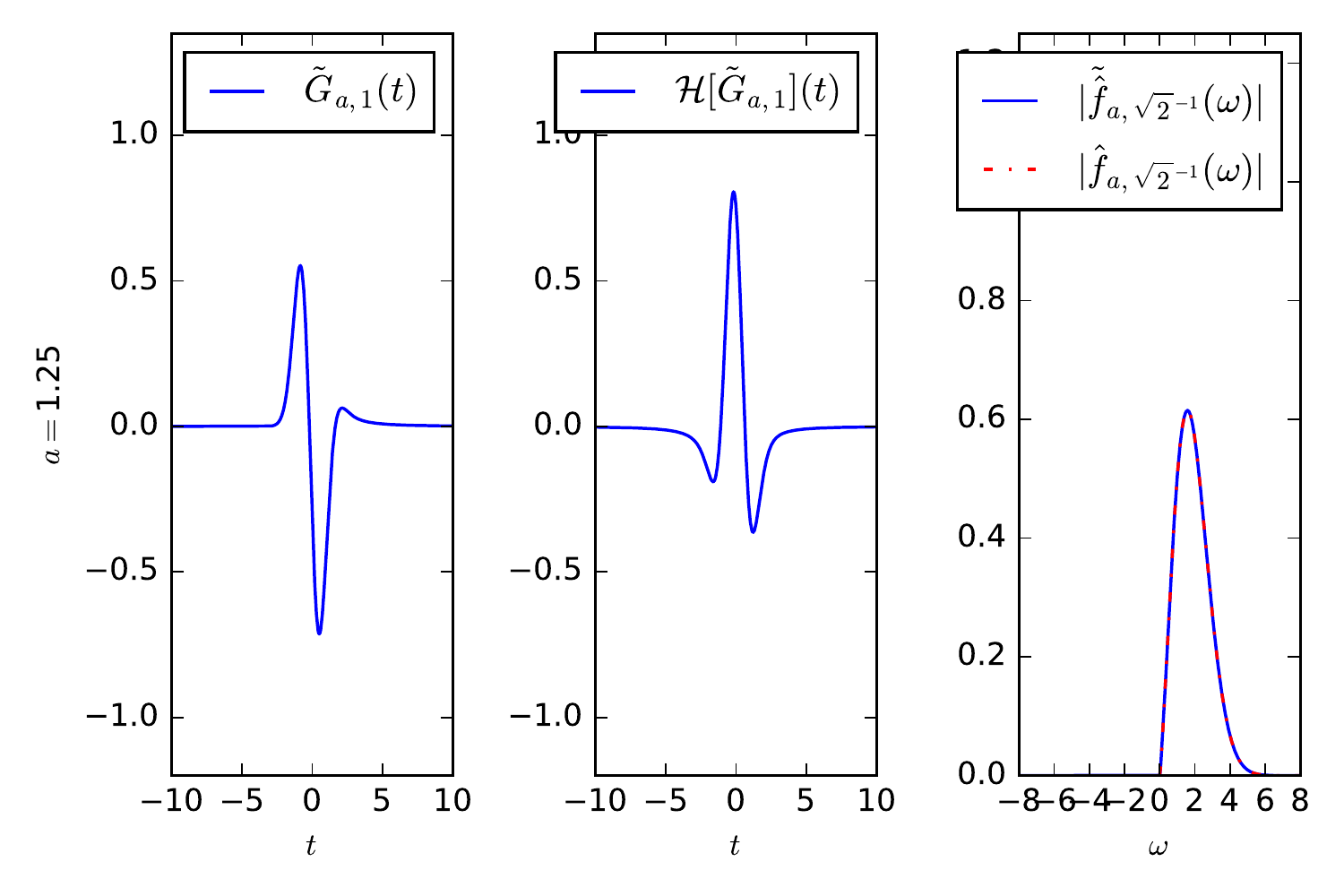} &
\includegraphics[width=0.3\textwidth]{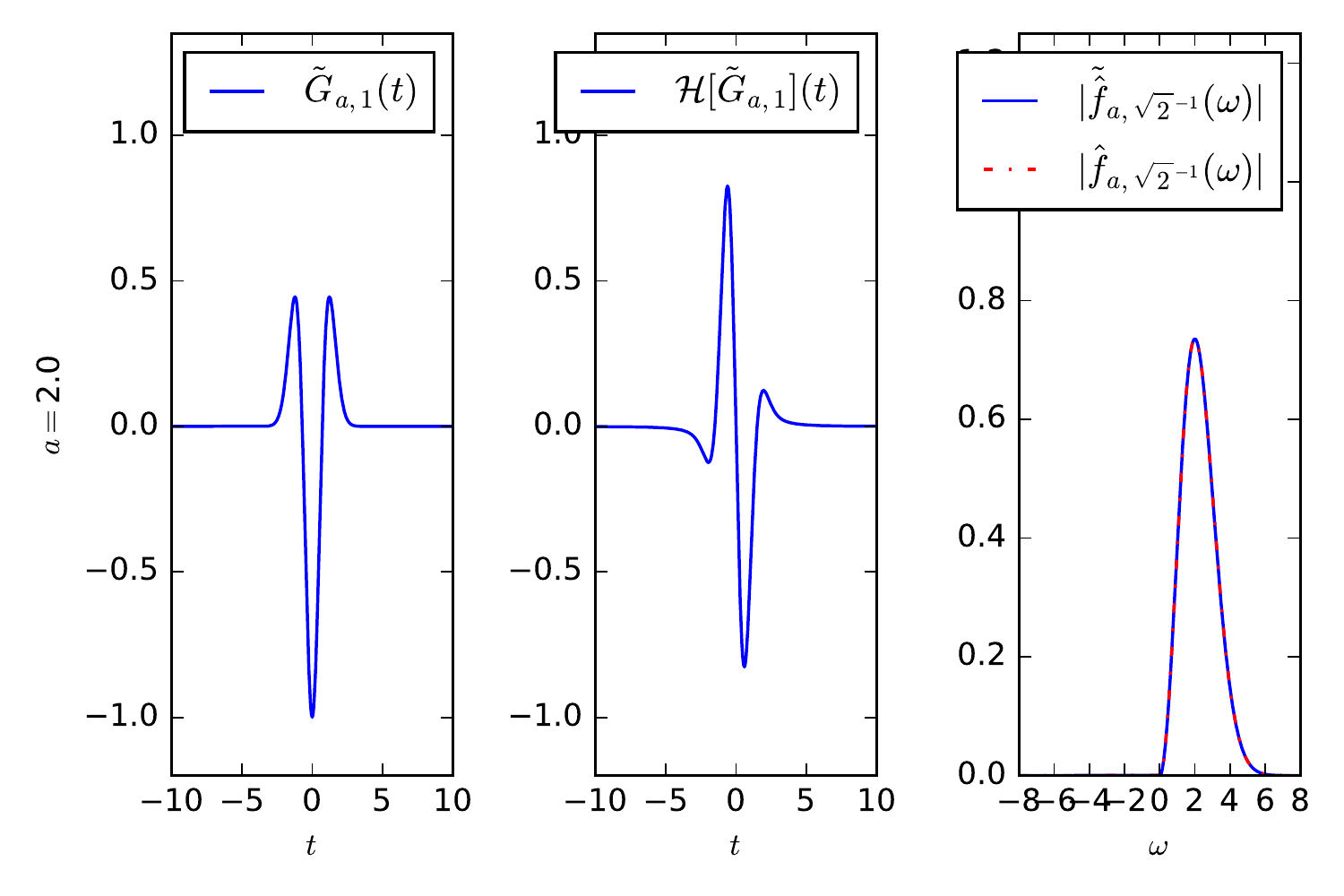}\tabularnewline
\hline 
\includegraphics[width=0.3\textwidth]{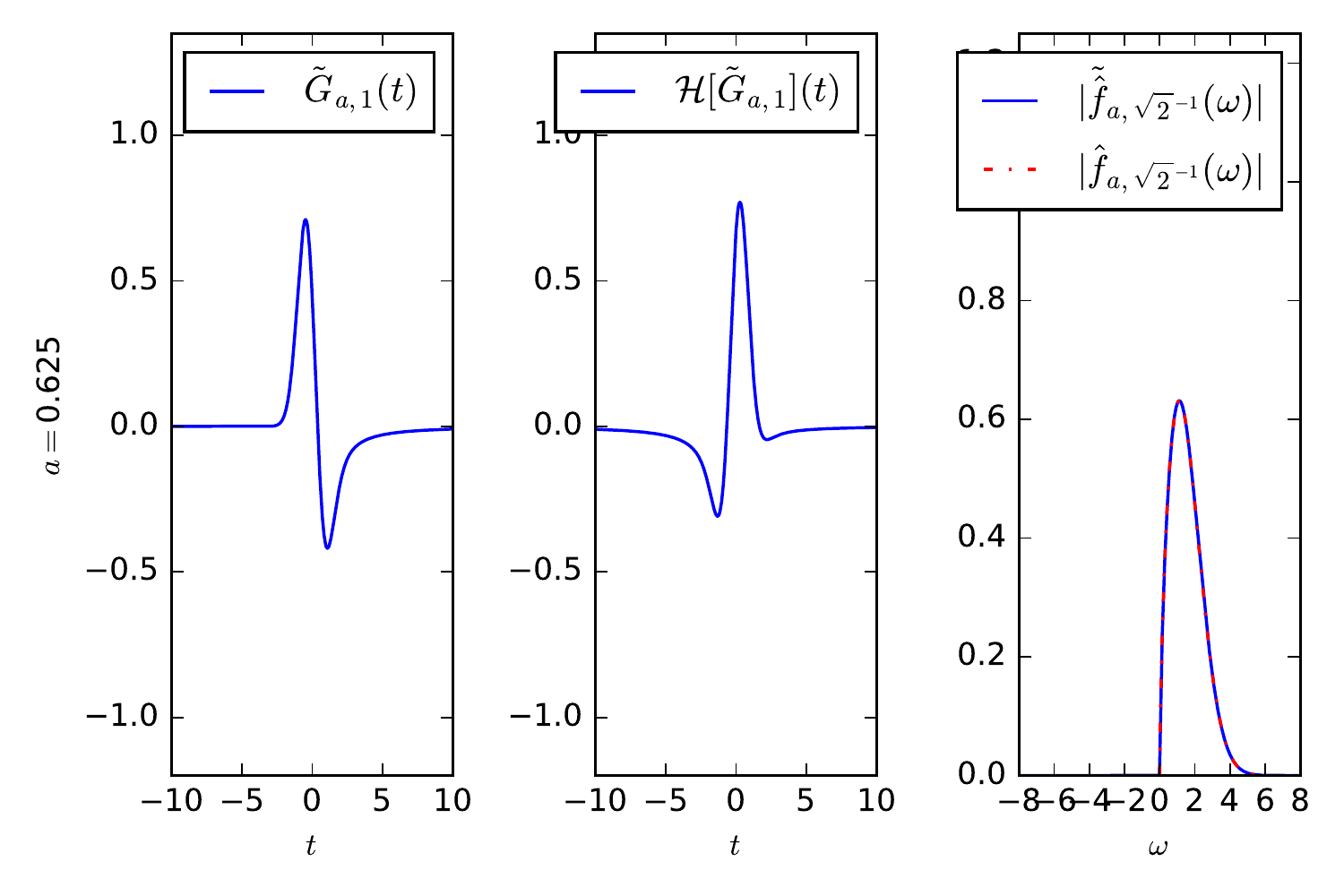} &
\includegraphics[width=0.3\textwidth]{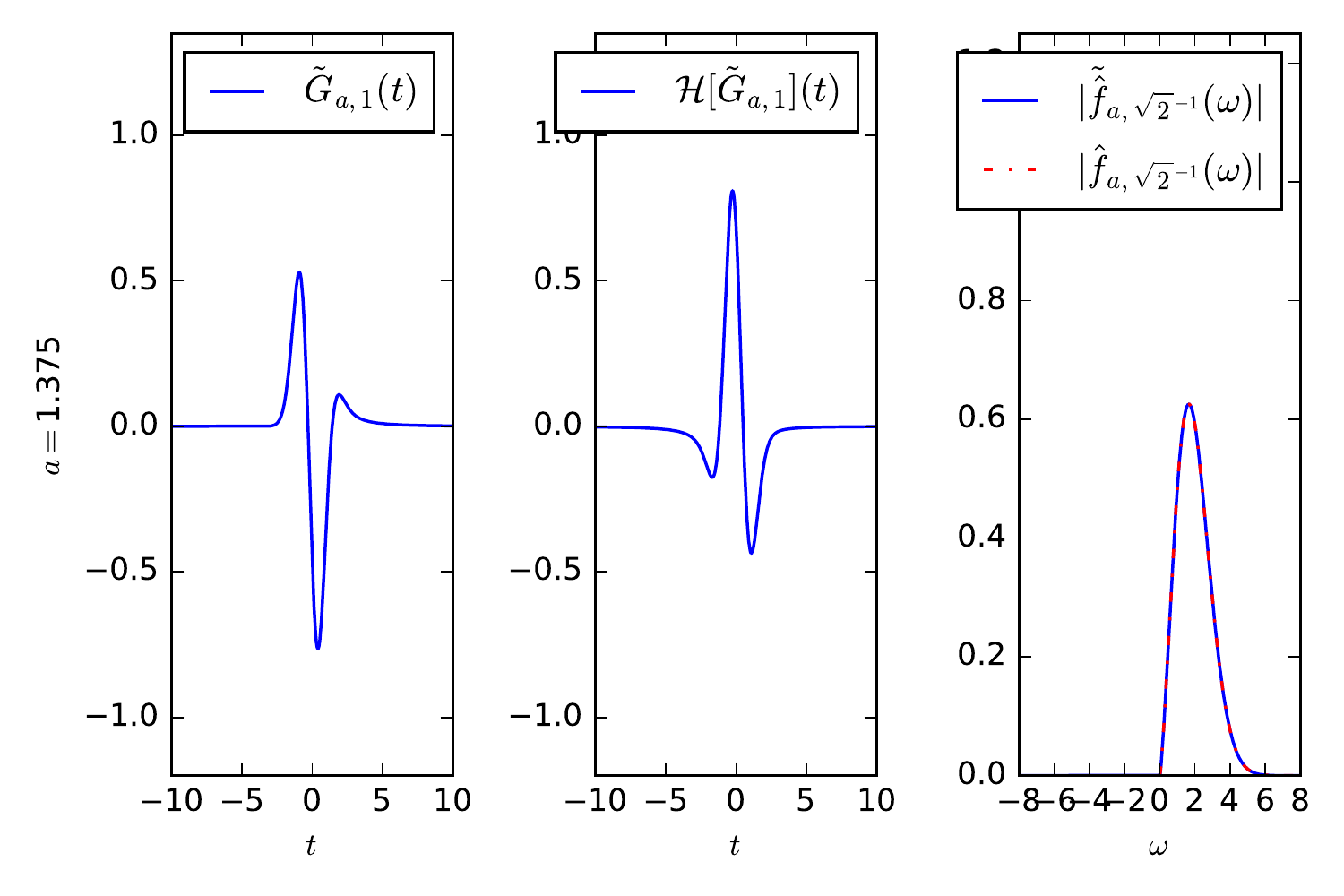} &
\includegraphics[width=0.3\textwidth]{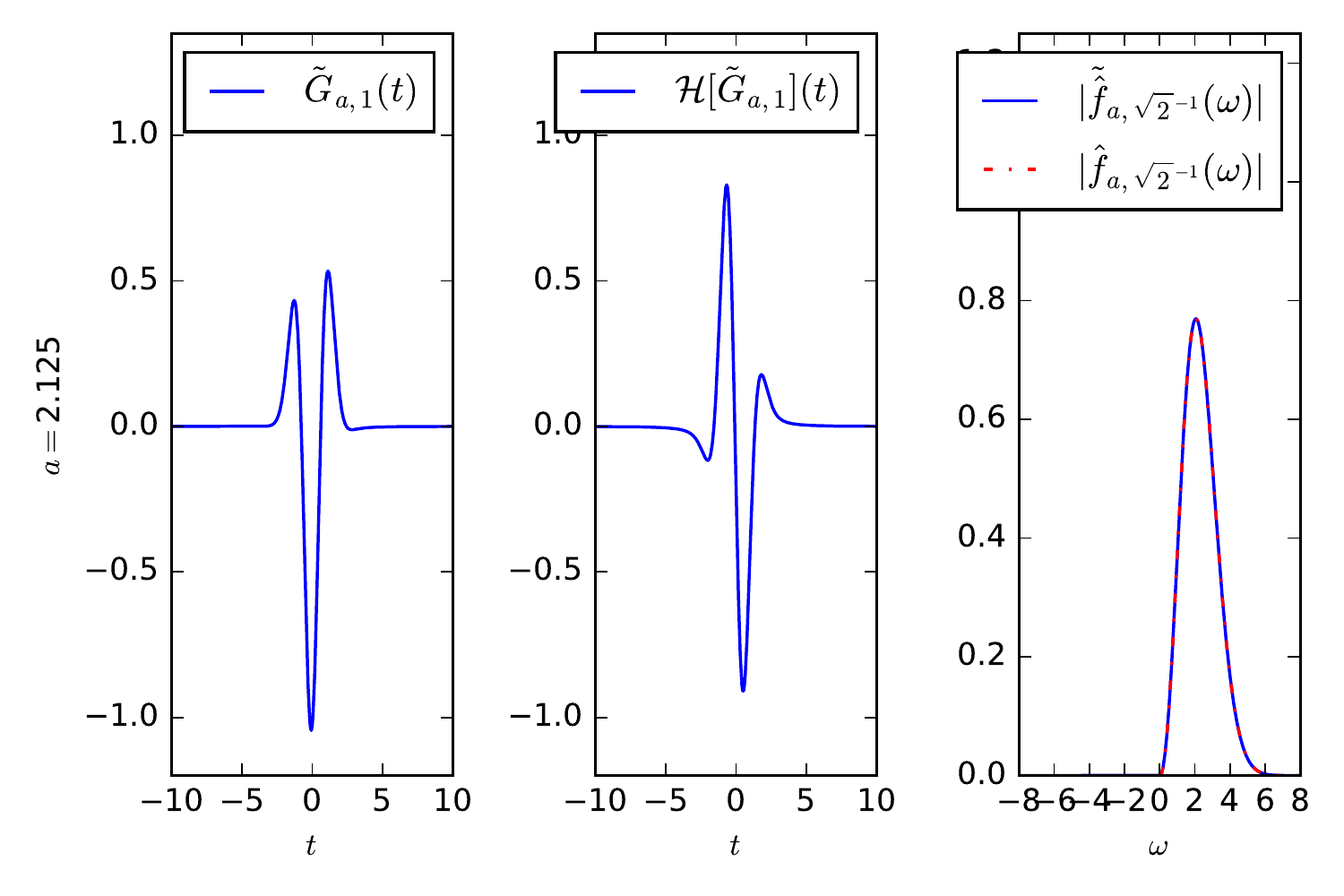}\tabularnewline
\hline 
\end{tabular}\caption{\label{fig:Plots-of-fractional-derivatives} Table of plots of $\tilde{G}_{a,\sigma}\left(t\right)$
(left), $\mathcal{H}\left[\tilde{G}_{a,\sigma}\right]\left(t\right)$
(middle), $\left|\tilde{\hat{f}}_{a,\sigma}\left(\omega\right)\right|$
(right, solid blue) and $\left|\hat{f}_{a,\sigma}\left(\omega\right)\right|$
(right, dashed-dotted  red)  for $\sigma=\sqrt{2}^{-1}$ and various
$a$ starting from $0$ to $2.125$ with increments of $0.125$. Along
each column of the table, $a$ increases top to bottom by $0.125$.
Along each row of the table, $a$ increases left to right by $0.75$. }
\end{figure}

We will denote the approximation of the fractional derivative of $f_{0,\sigma}\left(t\right)$
obtained through solution of the moment problem (\ref{eq:frac_derivative-moment_problem})
by 
\begin{align*}
\tilde{f}_{a,\sigma}\left(t\right) & =\sum_{m=1}^{M}\alpha_{m}g\left(\gamma_{m}t\right).
\end{align*}
{\small{} $\tilde{f}_{a,\sigma}\left(t\right)$} can also be written
as 
\begin{align*}
\tilde{f}_{a,\sigma}\left(t\right) & =\sqrt{\frac{2}{\pi}}\int_{0}^{\infty}\tilde{\hat{f}}_{a,\sigma}\left(\omega\right)\exp\left(\mathrm{i}\frac{\omega}{\sigma}t\right)d\omega=\tilde{G}_{a,\sigma}\left(t\right)+\mathrm{i}\mathcal{H}\left[\tilde{G}_{a,\sigma}\right]\left(t\right),
\end{align*}
for some $\tilde{\hat{f}}_{a,\sigma}\left(\omega\right)$ and $\tilde{G}_{a,\sigma}\left(t\right)=\mathrm{Re}\left\{ \tilde{f}_{a,\sigma}\left(t\right)\right\} $,
$\mathcal{H}\left[\tilde{G}_{a,\sigma}\right]\left(t\right)=\mathrm{Im}\left\{ \tilde{f}_{a,\sigma}\left(t\right)\right\} $.
Comparing with (\ref{eq:f_a_sigma_fourier}), $\tilde{\hat{f}}_{a,\sigma}\left(\omega\right)$
approximates $\hat{f}_{a,\sigma}\left(\omega\right)=\mathrm{e}^{-\frac{\omega^{2}}{2}}\left(\frac{\omega}{\sigma}\right)^{a}\exp\left(\mathrm{i}\,a\frac{\pi}{2}\right)$,
$\omega\ge0$. We present plots of $\tilde{G}_{a,\sigma}\left(t\right)$,
$\mathcal{H}\left[\tilde{G}_{a,\sigma}\right]\left(t\right)$ , $\left|\tilde{\hat{f}}_{a,\sigma}\left(\omega\right)\right|$
and $\left|\hat{f}_{a,\sigma}\left(\omega\right)\right|$ for various
values of $a$ in Figure \ref{fig:Plots-of-fractional-derivatives}.
In our computations, we used a rational approximation of the Dawson's
integral which is provided in Appendix \ref{sec:A-rational-approximation2dawson}.

Our approach can be generalized to approximate partial derivative
with respect to the order of fractional derivative which will be discussed
in Section \ref{sec:Partial-a-of-frac_der_a}.

\begin{table}
\begin{centering}
\begin{tabular}{|c|c|c|}
\hline 
$m$  &
$\alpha_{m}$  &
$\gamma_{m}$\tabularnewline
\hline 
\hline 
{\footnotesize{}1 } &
{\footnotesize{}-0.002327462216272 - 0.002323117038281i} &
{\footnotesize{}0.022707194026268 - 0.000014485869526i}\tabularnewline
\hline 
{\footnotesize{}2 } &
{\footnotesize{}-0.008686498484022 - 0.008672831744140i} &
{\footnotesize{}0.088058420558651 - 0.000050483776843i}\tabularnewline
\hline 
{\footnotesize{}3 } &
{\footnotesize{}-0.017746649906246 - 0.017725233118614i} &
{\footnotesize{}0.188199334944285 - 0.000090762508709i}\tabularnewline
\hline 
{\footnotesize{}4 } &
{\footnotesize{}-0.028479502359746 - 0.028455359343906i} &
{\footnotesize{}0.311644959787523 - 0.000118301948941i}\tabularnewline
\hline 
{\footnotesize{}5 } &
{\footnotesize{}-0.040781085680949 - 0.040758739735612i} &
{\footnotesize{}0.445550822358306 - 0.000124204777278i}\tabularnewline
\hline 
{\footnotesize{}6 } &
{\footnotesize{}-0.055691041043028 - 0.055672403045923i} &
{\footnotesize{}0.578033122632225 - 0.000110160892971i}\tabularnewline
\hline 
{\footnotesize{}7 } &
{\footnotesize{}-0.075799236577609 - 0.075783432858215i} &
{\footnotesize{}0.699719488265476 - 0.000084621403249i}\tabularnewline
\hline 
{\footnotesize{}8 } &
{\footnotesize{}-0.106825493665143 - 0.106809444597271i} &
{\footnotesize{}0.804243109619030 - 0.000056714249345i}\tabularnewline
\hline 
{\footnotesize{}9 } &
{\footnotesize{}2.224002927119579 + 2.223687081665376i} &
{\footnotesize{}1.001024232778366 + 0.000000290415923i}\tabularnewline
\hline 
{\footnotesize{}10 } &
{\footnotesize{}-0.846706907066594 - 0.846588952992167i} &
{\footnotesize{}0.986812368161672 - 0.000003755178001i}\tabularnewline
\hline 
{\footnotesize{}11 } &
{\footnotesize{}-0.296187315352343 - 0.296147582995713i} &
{\footnotesize{}0.948969382632748 - 0.000014668011082i}\tabularnewline
\hline 
{\footnotesize{}12 } &
{\footnotesize{}-0.163403417748510 - 0.163381667176418i} &
{\footnotesize{}0.887911302390365 - 0.000032513183747i}\tabularnewline
\hline 
\end{tabular}
\par\end{centering}
\caption{\label{tab:(alpha_m,gamma_m)_frac_gaussian_in_gaussian+dawson} Table
of $(\alpha_{m},\gamma_{m})$ obtained by solving the moment problem
(\ref{eq:frac_derivative-moment_problem}) for $a=2^{-1}$ and $\sigma=\sqrt{2}^{-1}$.
Approximations $\tilde{f}_{2^{-1},\sqrt{2}^{-1}}\left(t\right)$ and
$\tilde{f}_{a,\sqrt{2}^{-1}}\left(t\right)$ for other values of$a$
are presented in Figure \ref{fig:Plots-of-fractional-derivatives}}
\end{table}

\subsection{Approximation error}

Recall that 
\begin{align*}
g\left(t\right) & =f_{0,\sigma}\left(\sqrt{2}t\right)=\sqrt{\frac{2}{\pi}}\int_{0}^{\infty}\mathrm{e}^{-\frac{\omega^{2}}{2}}\exp\left(\mathrm{i}\frac{\omega}{\sqrt{2}}t\right)d\omega
\end{align*}
 and 
\begin{align*}
f_{a,\sigma}\left(t\right) & =\sqrt{\frac{2}{\pi}}\int_{0}^{\infty}\mathrm{e}^{-\frac{\omega^{2}}{2}}\left(\frac{\omega}{\sigma}\right)^{a}\exp\left(\mathrm{i}\left[\frac{\omega}{\sigma}t+a\frac{\pi}{2}\right]\right)d\omega.
\end{align*}
 Then 
\begin{align*}
\tilde{f}_{a,\sigma}\left(t\right) & =\sum_{m=1}^{M}\alpha_{m}g\left(\gamma_{m}t\right)\\
 & =\sqrt{\frac{2}{\pi}}\int_{0}^{\infty}\mathrm{e}^{-\frac{\omega^{2}}{2}}\sum_{m=1}^{M}\alpha_{m}\exp\left(\mathrm{i}\frac{\gamma_{m}}{\sqrt{2}}\omega t\right)d\omega.
\end{align*}

Consequently, the error and its L2 norm are given by 
\begin{multline*}
\left[\tilde{f}_{a,\sigma}-f_{a,\sigma}\right]\left(t\right)\\
=\sqrt{\frac{2}{\pi}}\int_{0}^{\infty}\mathrm{e}^{-\frac{\omega^{2}}{2}}\left[\sum_{m=1}^{M}\alpha_{m}\exp\left(\mathrm{i}\frac{\gamma_{m}}{\sqrt{2}}\omega t\right)-\left(\frac{\omega}{\sigma}\right)^{a}\exp\left(\mathrm{i}\left[\frac{\omega}{\sigma}t+a\frac{\pi}{2}\right]\right)\right]d\omega
\end{multline*}
and 
\begin{multline*}
\int_{-\infty}^{\infty}\left|\left[\tilde{f}_{a,\sigma}-f_{a,\sigma}\right]\left(t\right)\right|^{2}dt\\
=4\left[\begin{array}{l}
\sum_{m'=1}^{M}\sum_{m=1}^{M}\alpha_{m'}^{*}\alpha_{m}\left[\frac{\gamma_{m}\gamma_{m'}^{*}}{2}\right]\sqrt{\frac{\pi}{\gamma_{m'}^{*2}+\gamma_{m}^{2}}}+\frac{\sigma^{1-2a}}{\pi}\Gamma\left(a+\frac{1}{2}\right)\\
-2\text{Re}\left\{ \frac{\sigma}{\sqrt{2}}\exp\left(-\mathrm{i}a\frac{\pi}{2}\right)\sum_{m=1}^{M}\alpha_{m}\gamma_{m}^{a}\left(1+\frac{\sigma^{2}\gamma_{m}^{2}}{2}\right)^{-\frac{a+1}{2}}\Gamma\left(\frac{a+1}{2}\right)\right\} 
\end{array}\right],
\end{multline*}
 for $a>-\frac{1}{2}$, respectively.

\section{Approximating derivative with respect to the order of the fractional
derivative\label{sec:Partial-a-of-frac_der_a}}

\begin{figure}
\begin{tabular}{|c|c|c|}
\hline 
\includegraphics[width=0.3\textwidth]{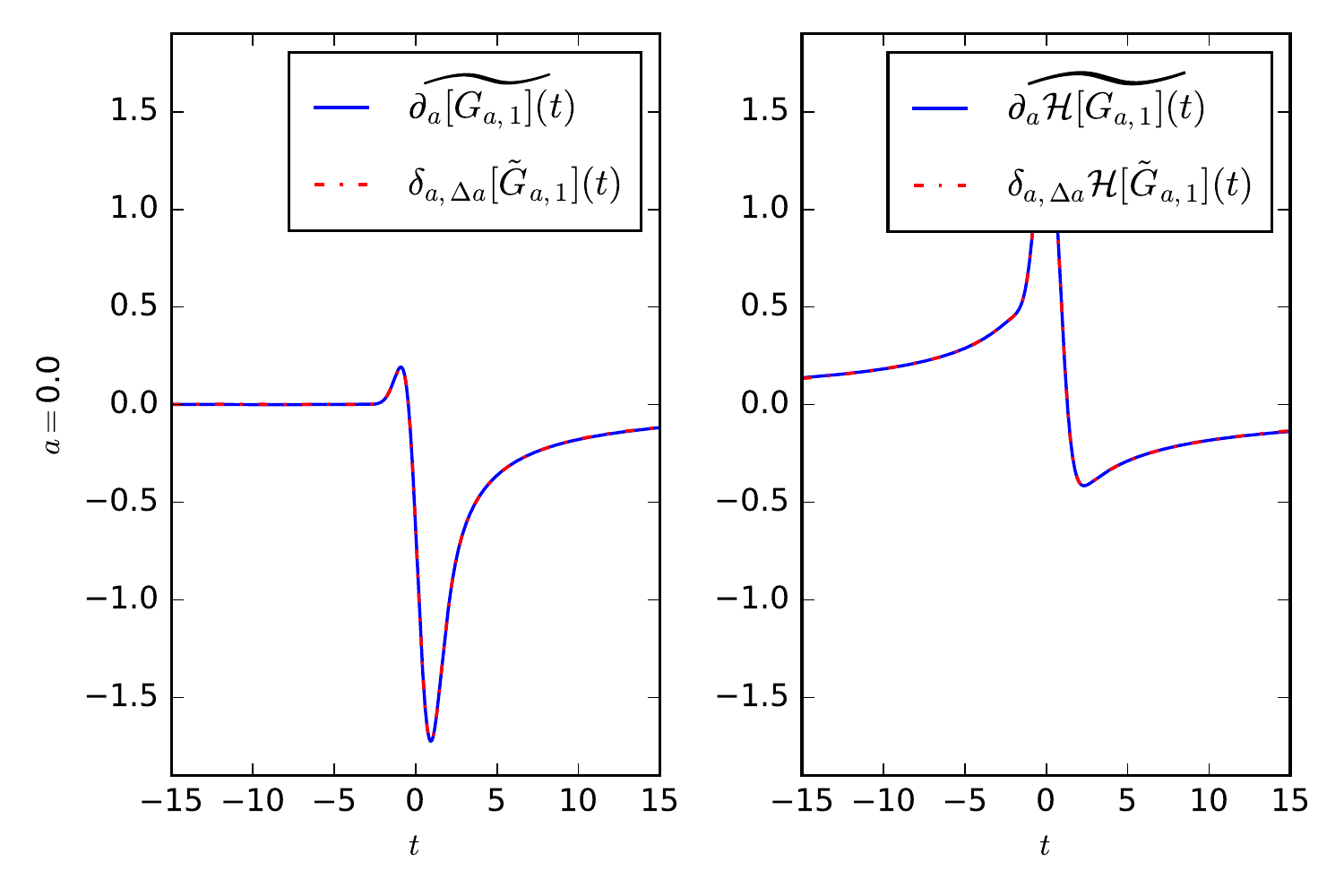} &
\includegraphics[width=0.3\textwidth]{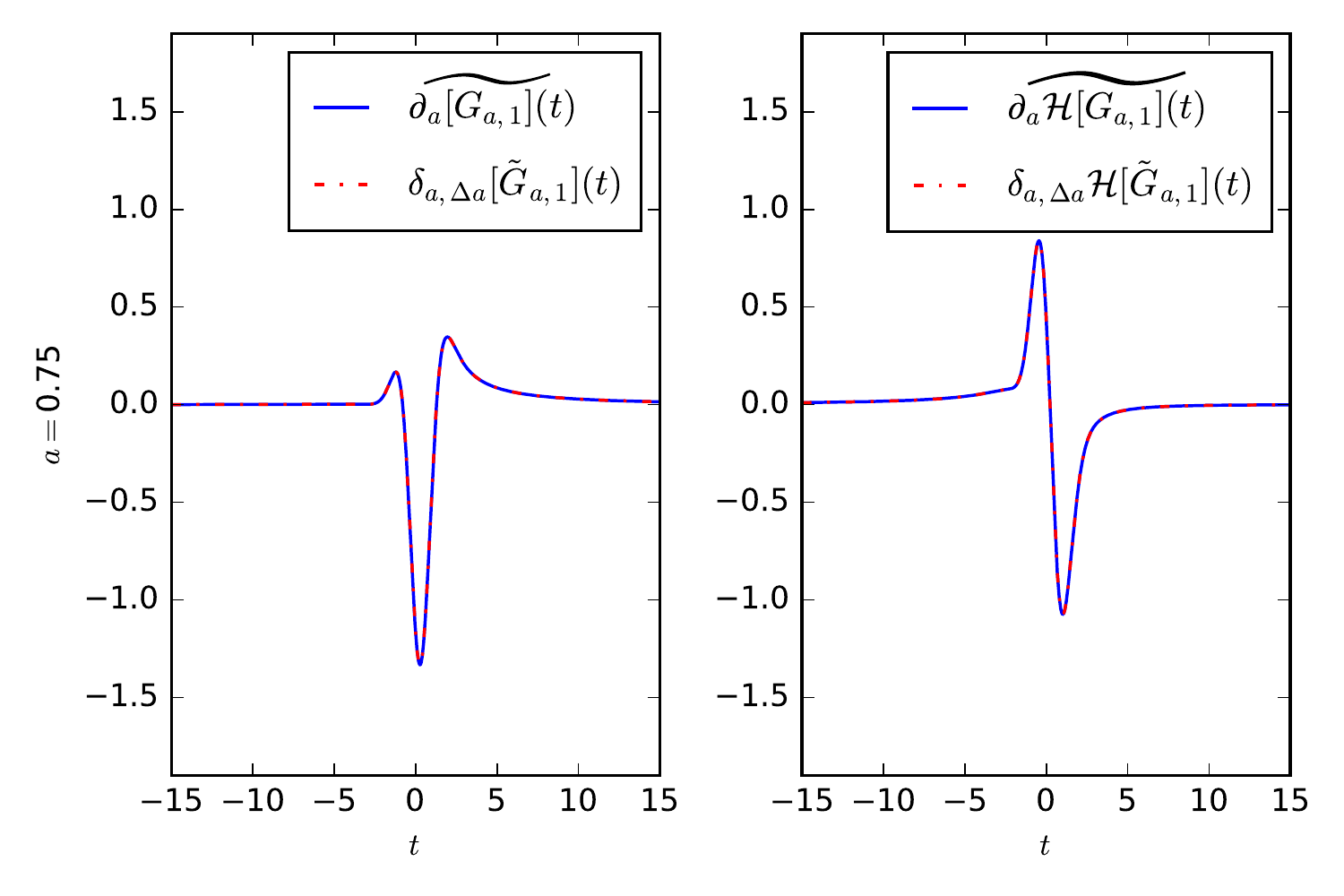} &
\includegraphics[width=0.3\textwidth]{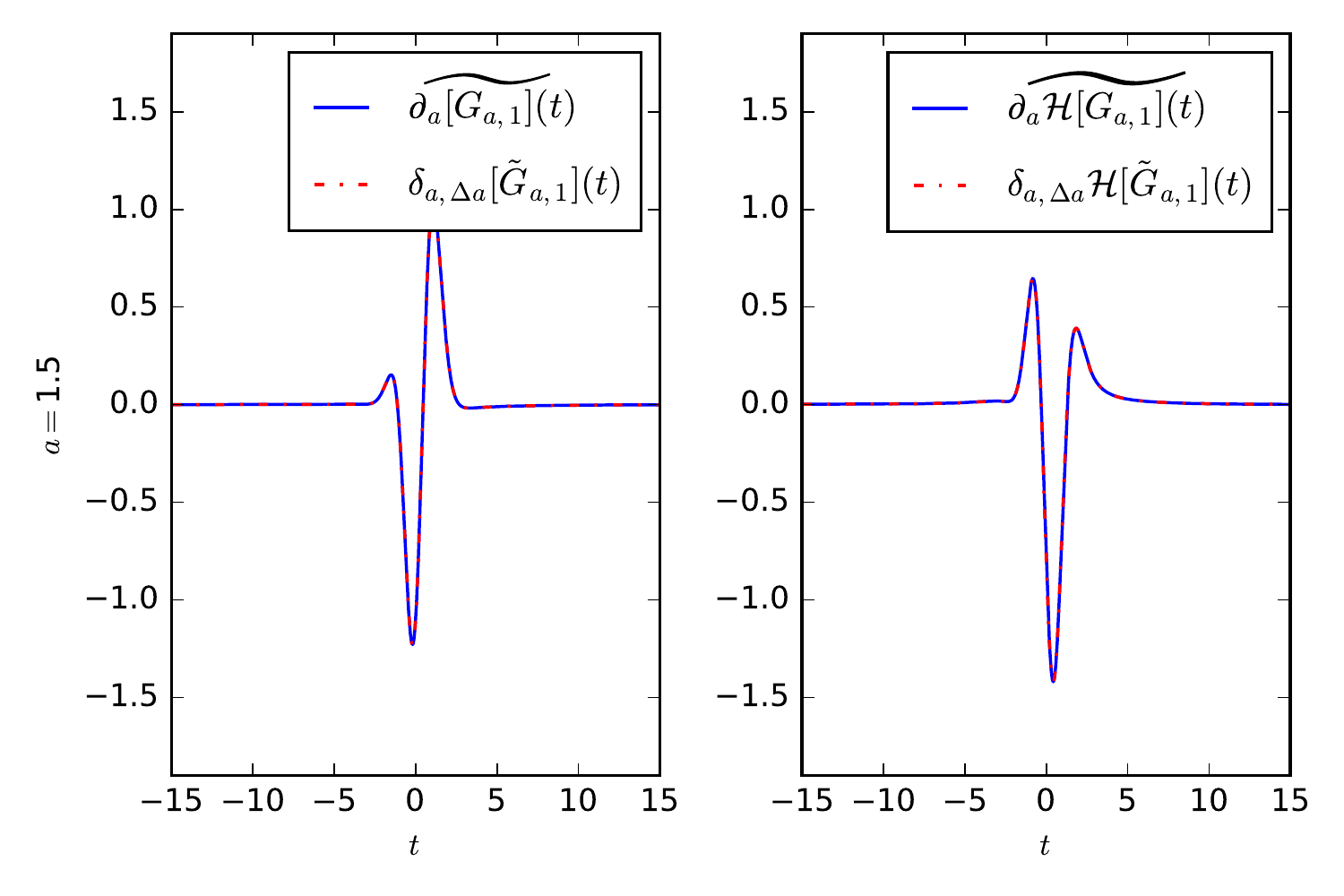}\tabularnewline
\hline 
\includegraphics[width=0.3\textwidth]{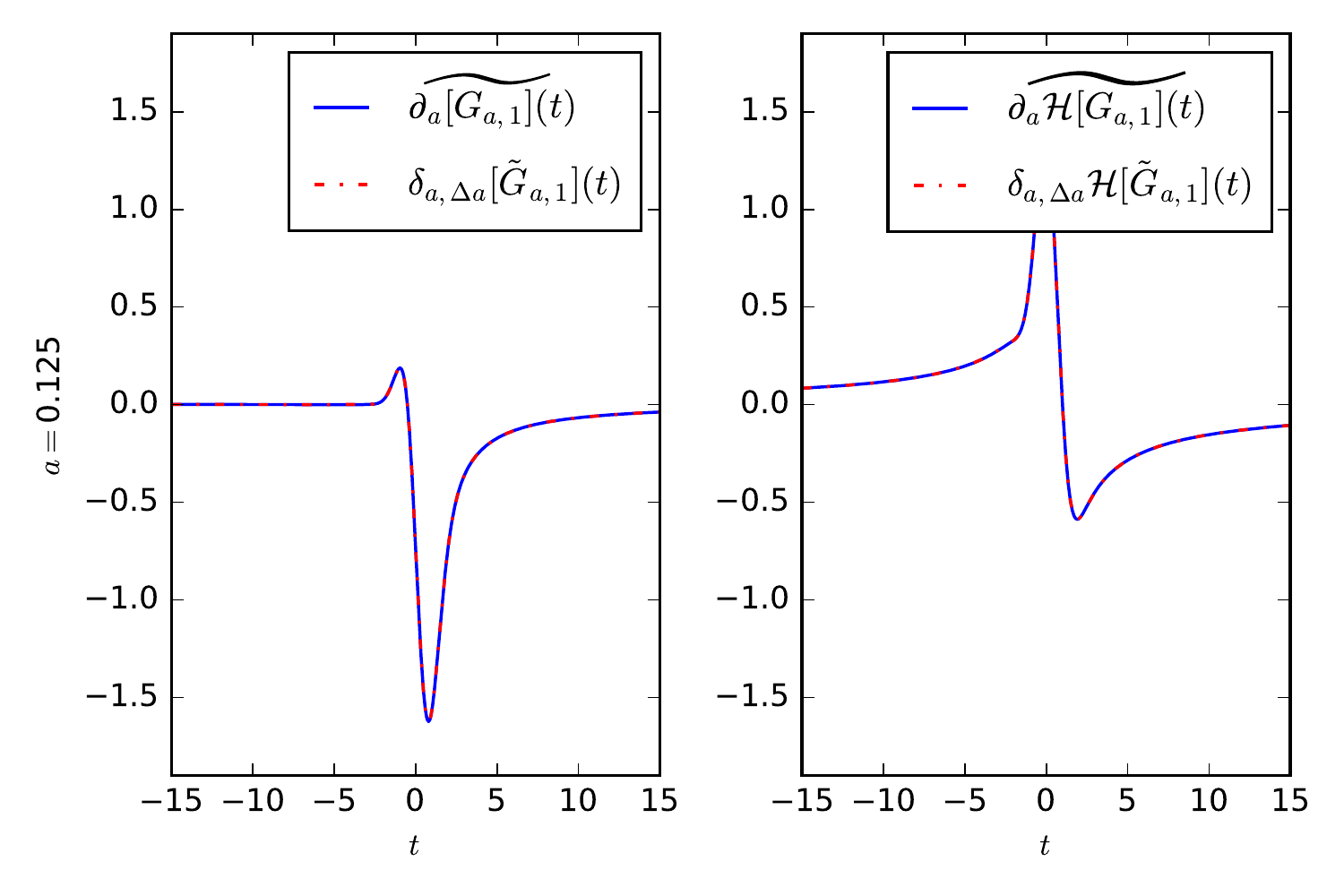} &
\includegraphics[width=0.3\textwidth]{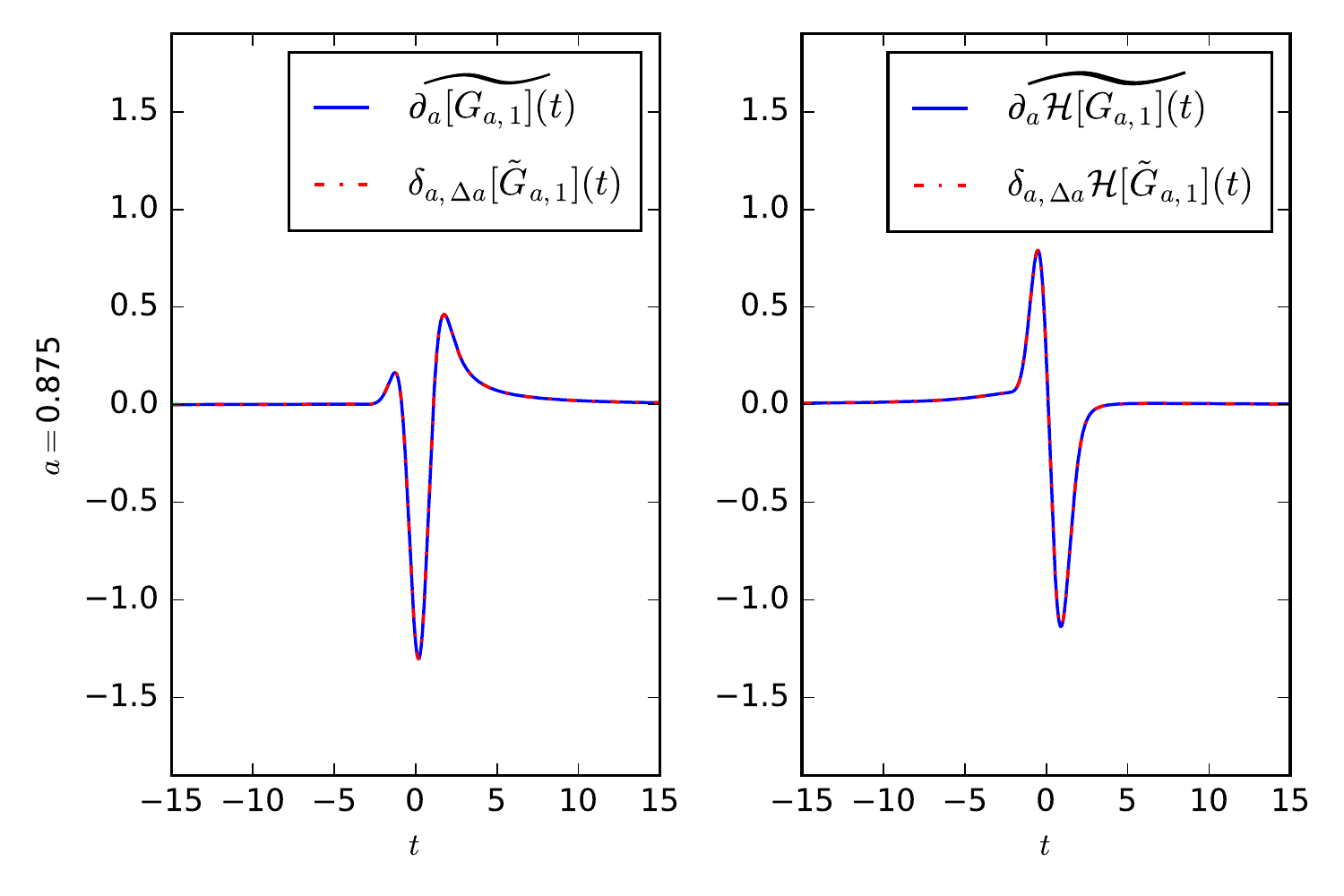} &
\includegraphics[width=0.3\textwidth]{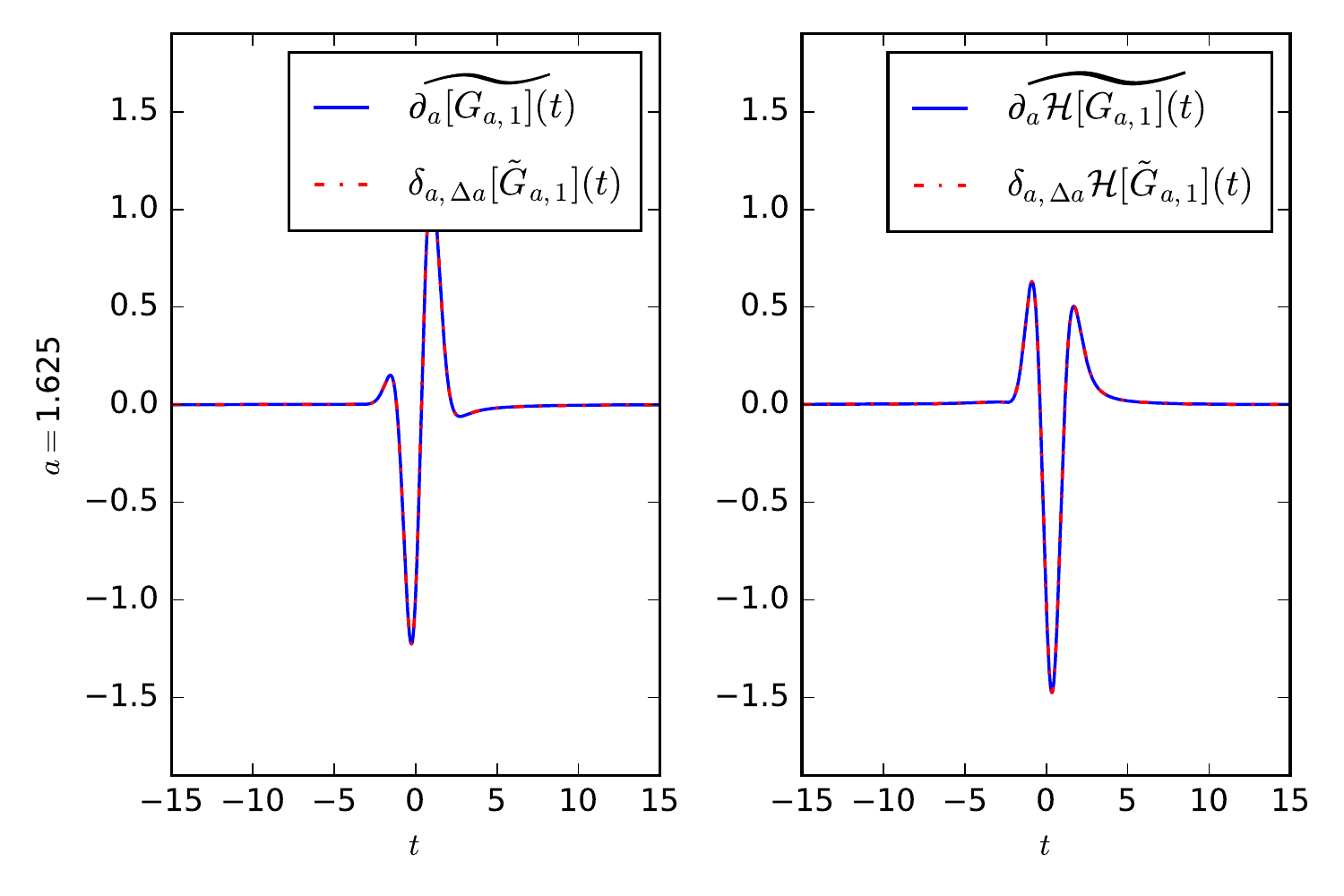}\tabularnewline
\hline 
\includegraphics[width=0.3\textwidth]{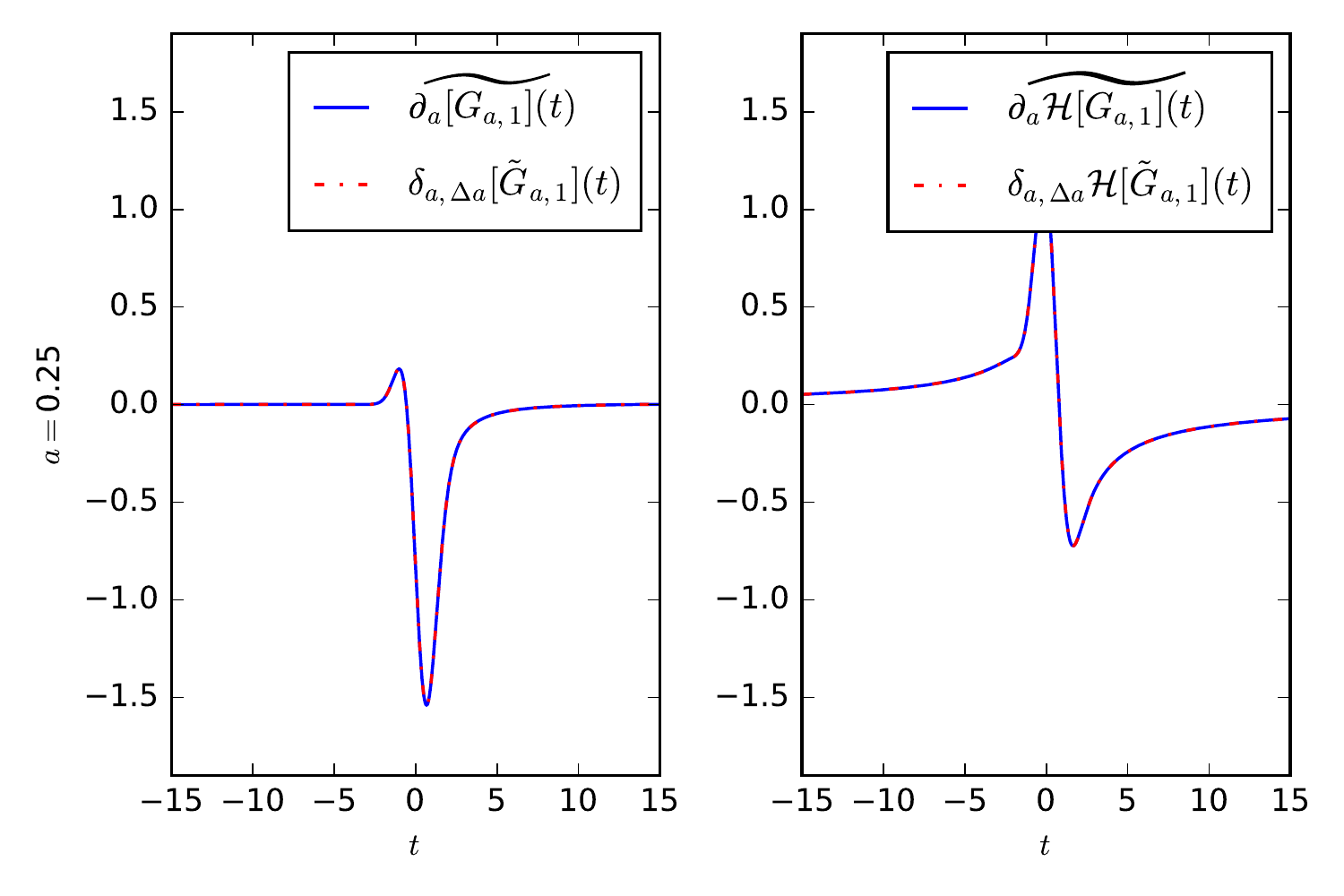} &
\includegraphics[width=0.3\textwidth]{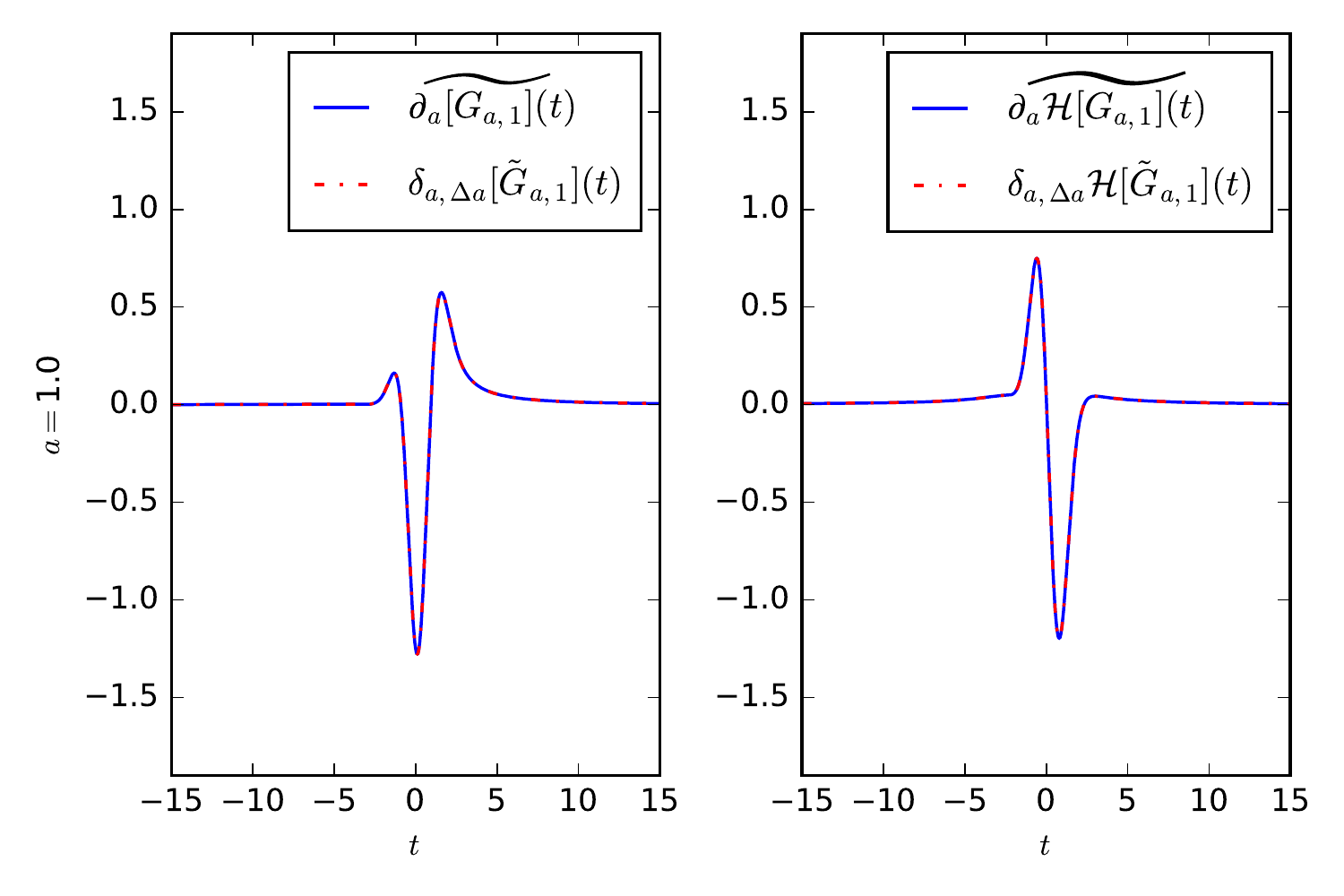} &
\includegraphics[width=0.3\textwidth]{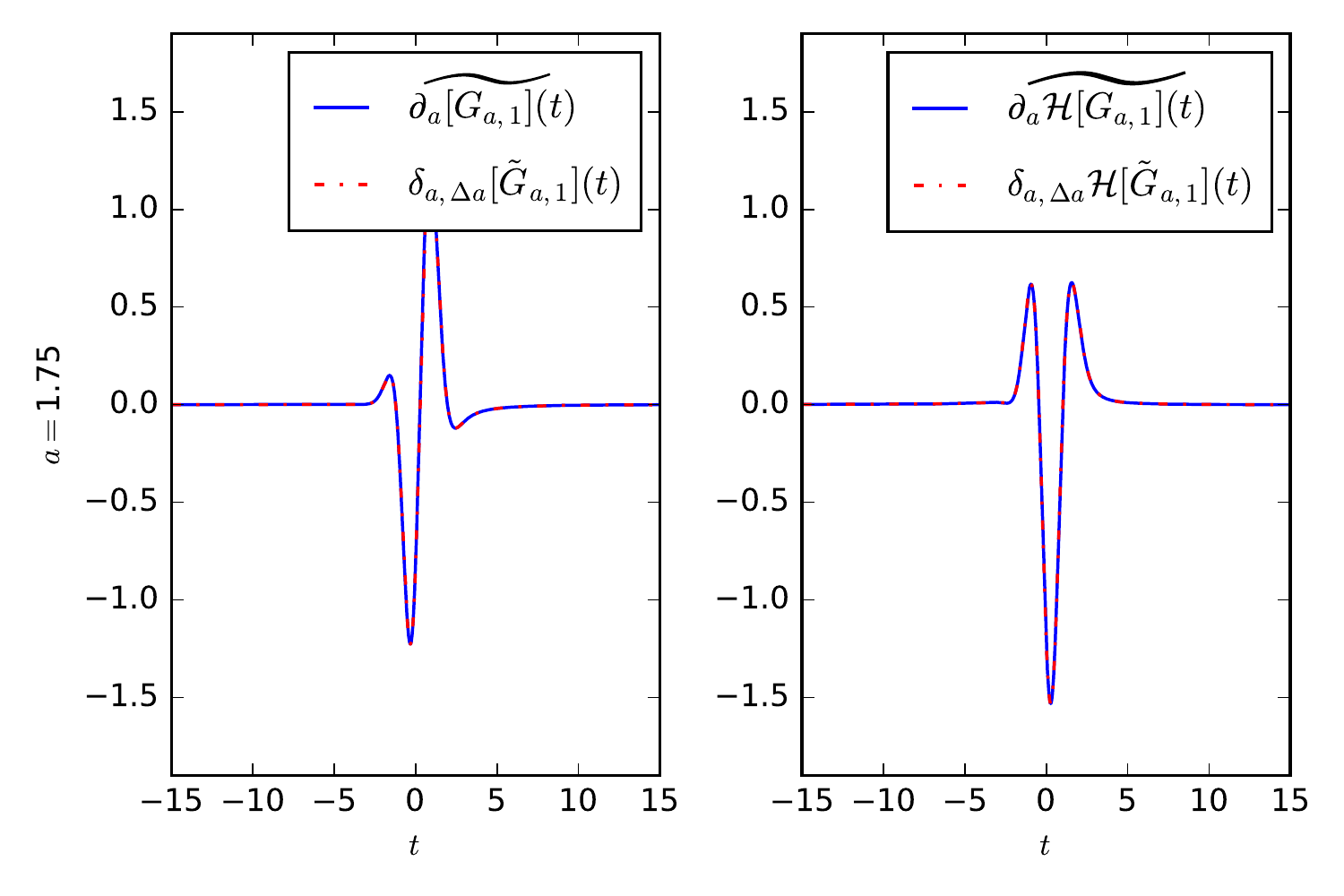}\tabularnewline
\hline 
\includegraphics[width=0.3\textwidth]{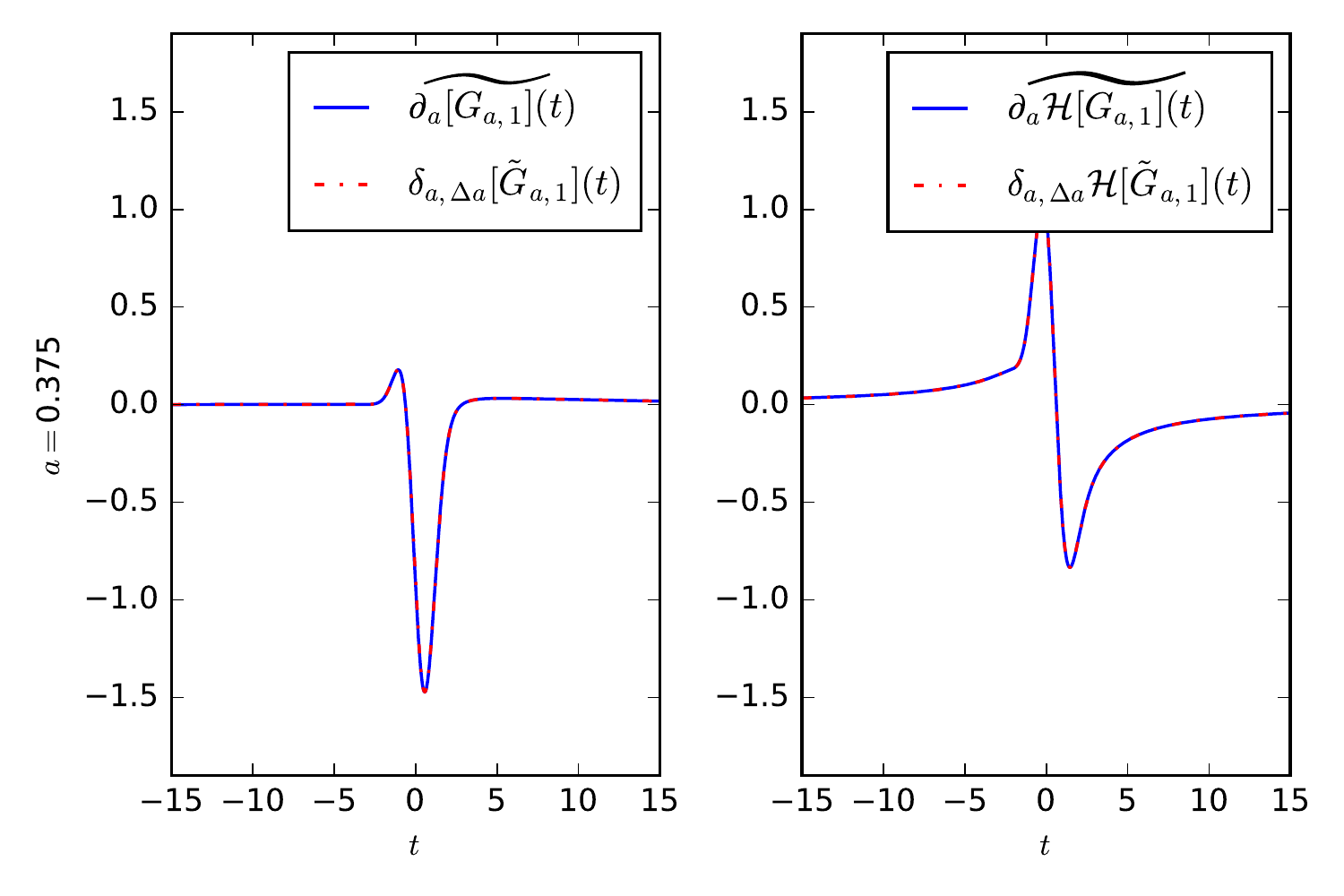} &
\includegraphics[width=0.3\textwidth]{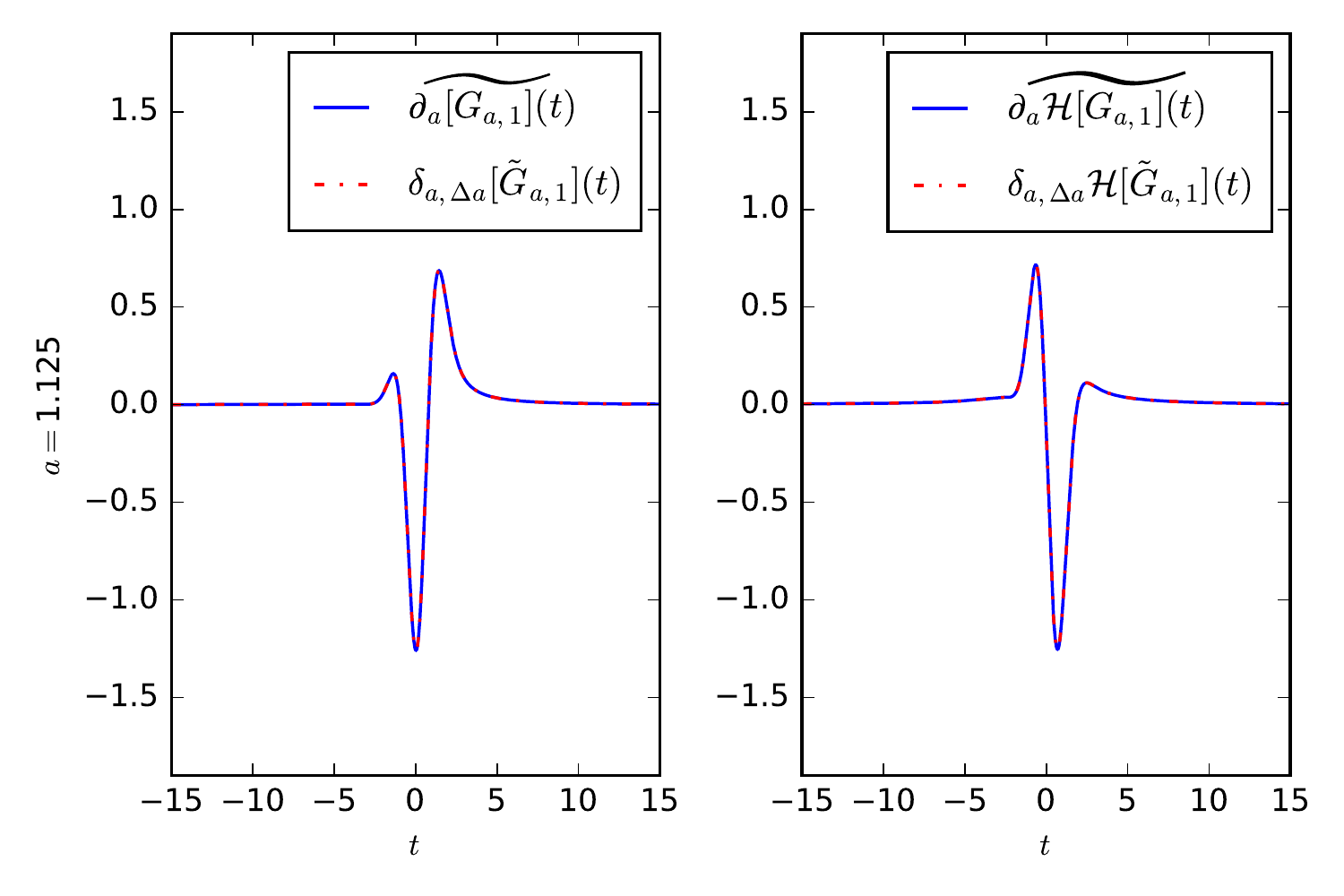} &
\includegraphics[width=0.3\textwidth]{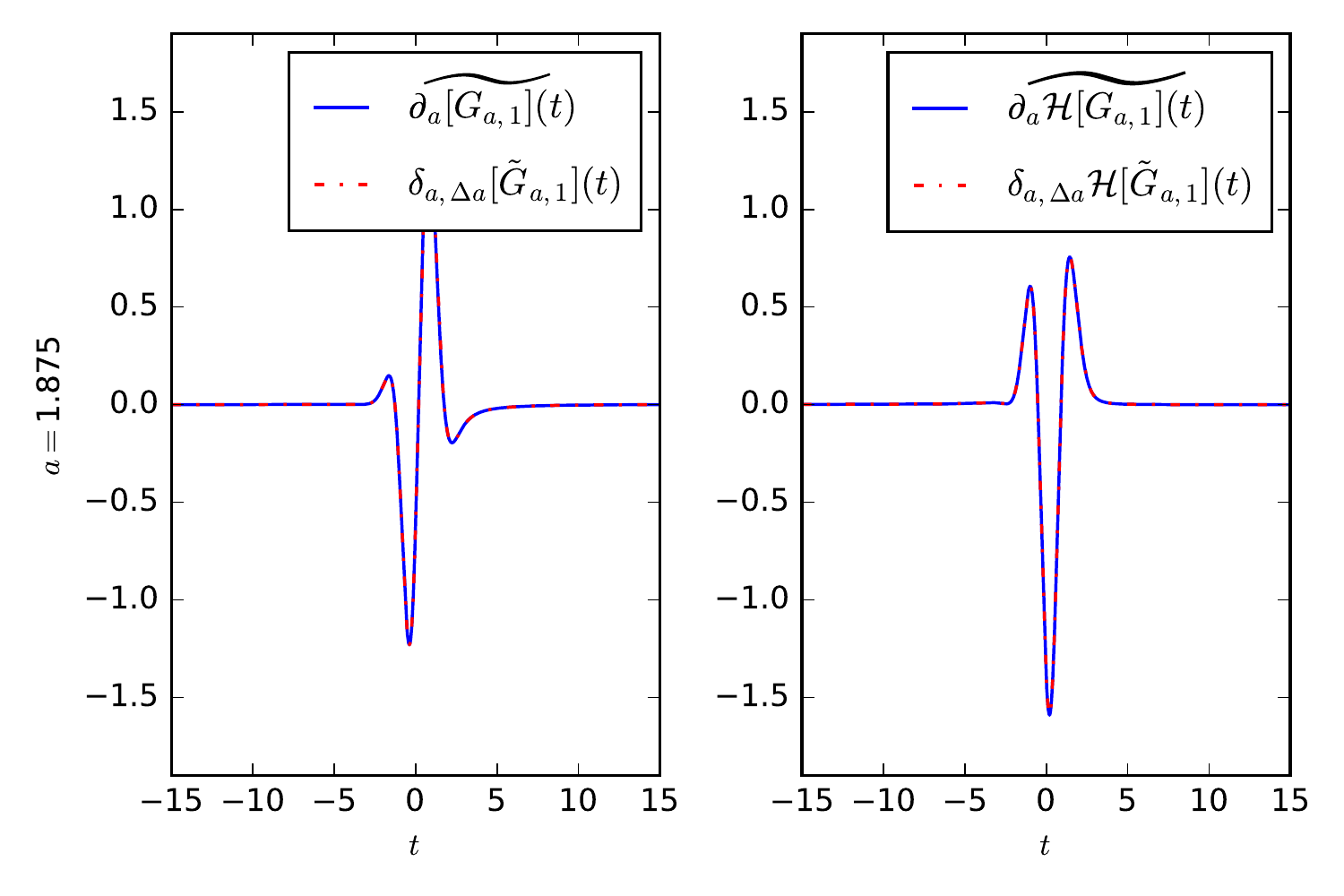}\tabularnewline
\hline 
\includegraphics[width=0.3\textwidth]{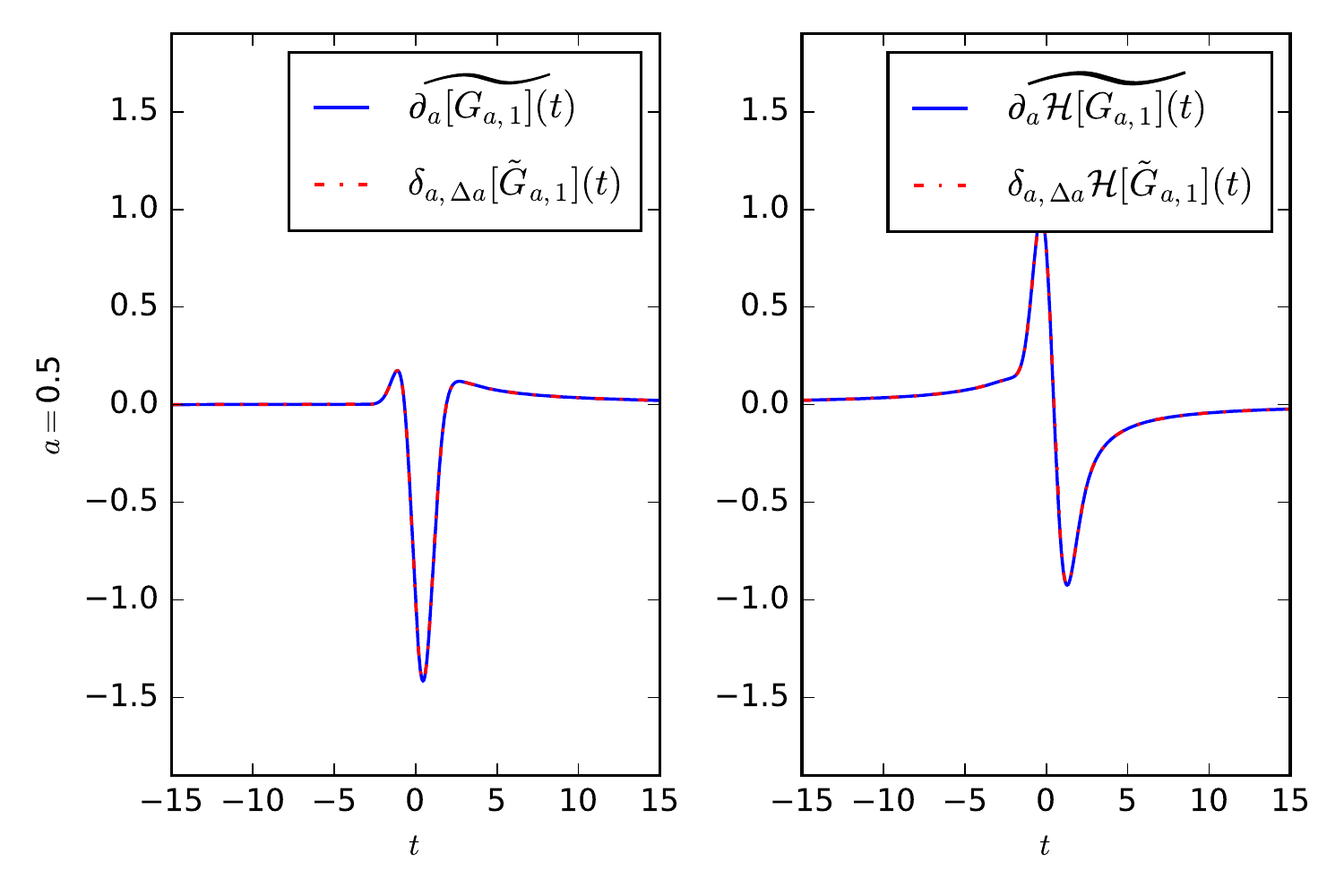} &
\includegraphics[width=0.3\textwidth]{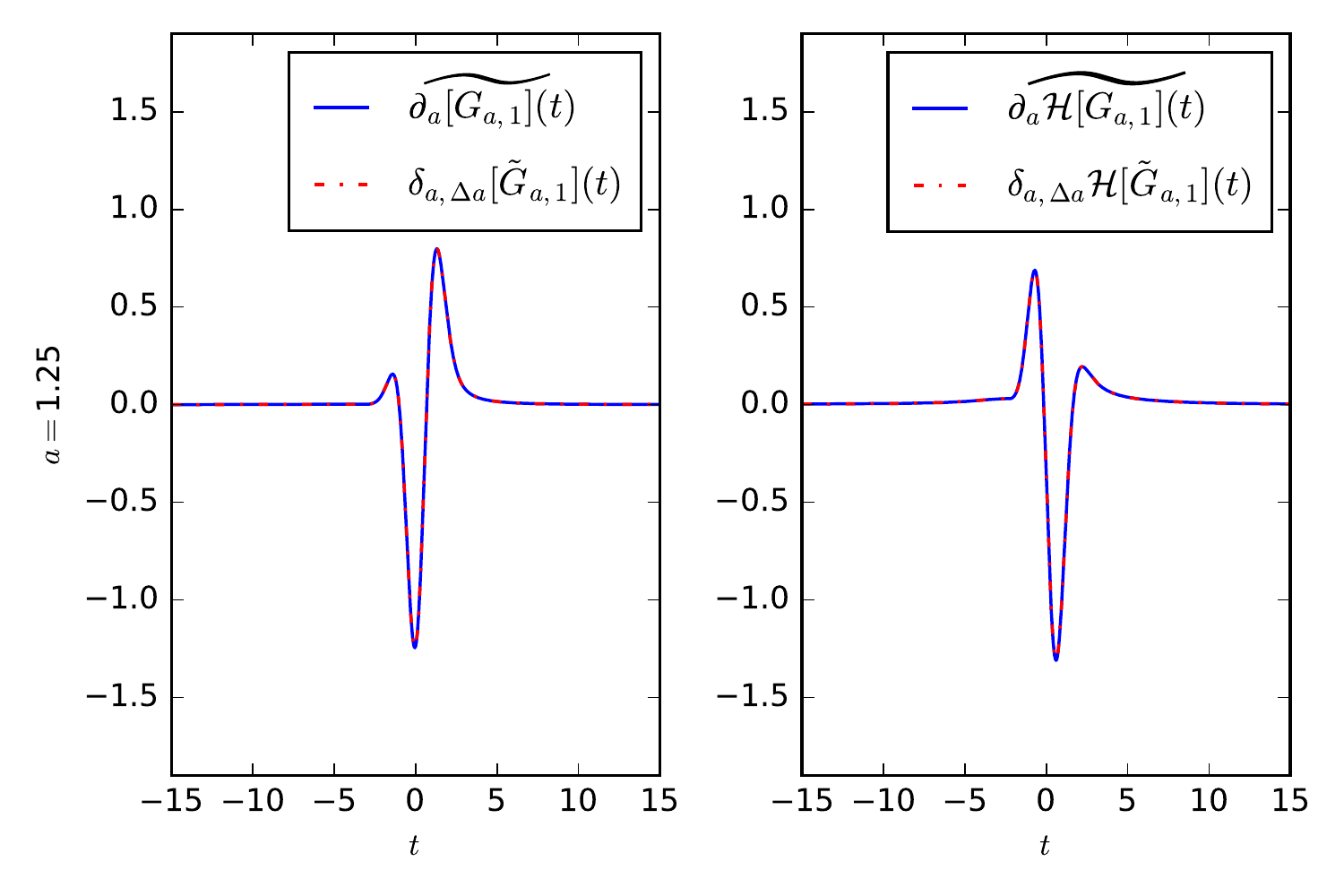} &
\includegraphics[width=0.3\textwidth]{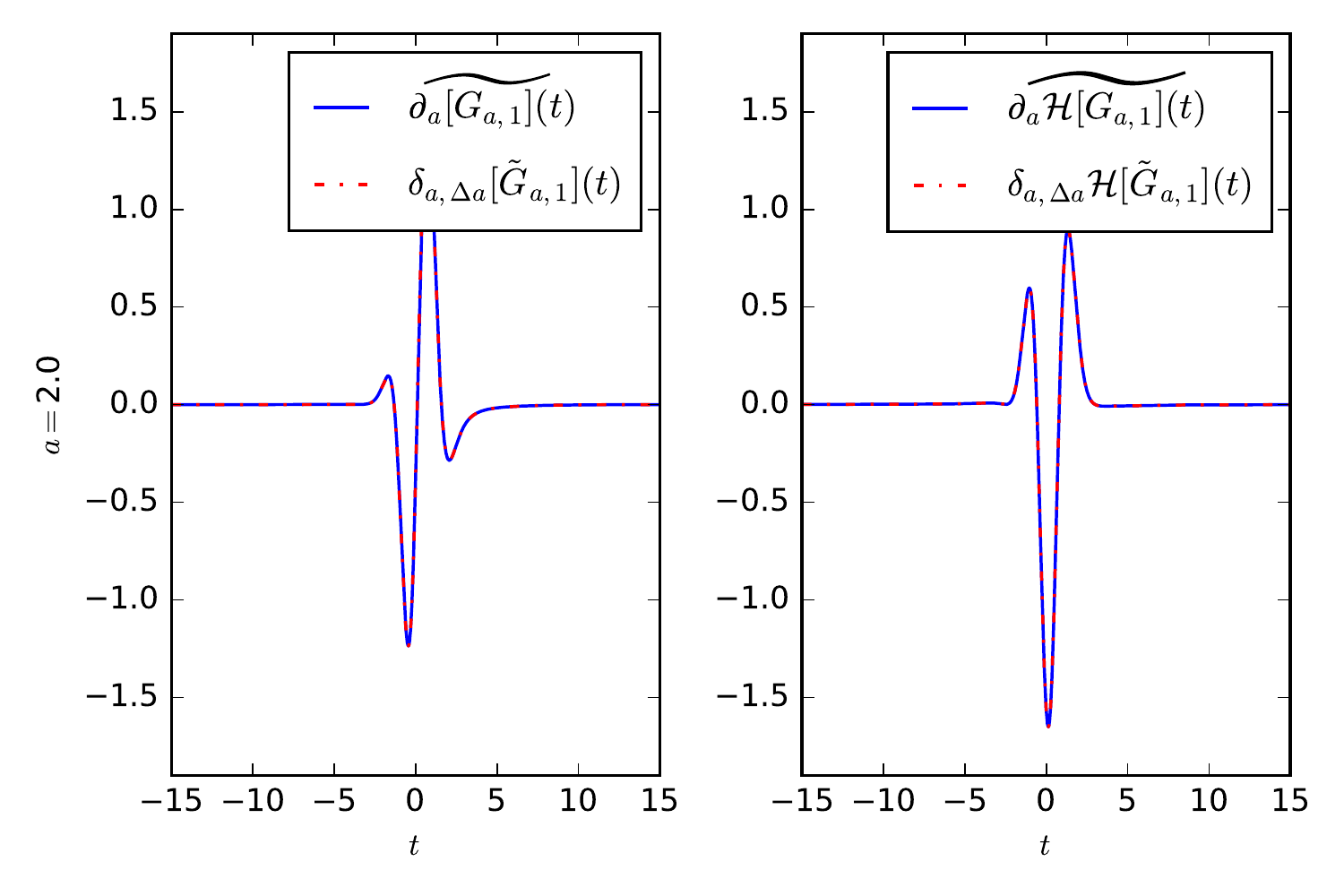}\tabularnewline
\hline 
\includegraphics[width=0.3\textwidth]{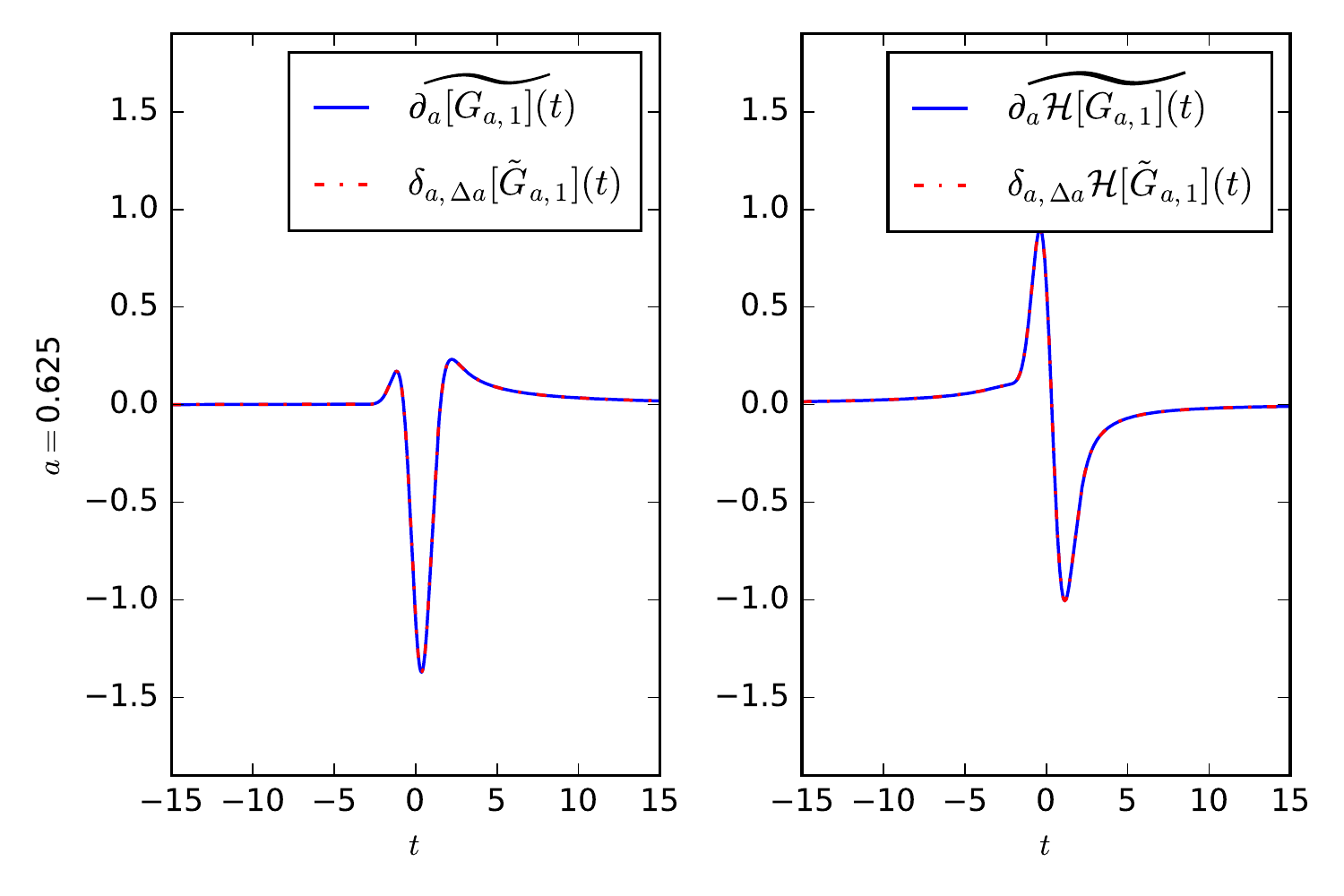} &
\includegraphics[width=0.3\textwidth]{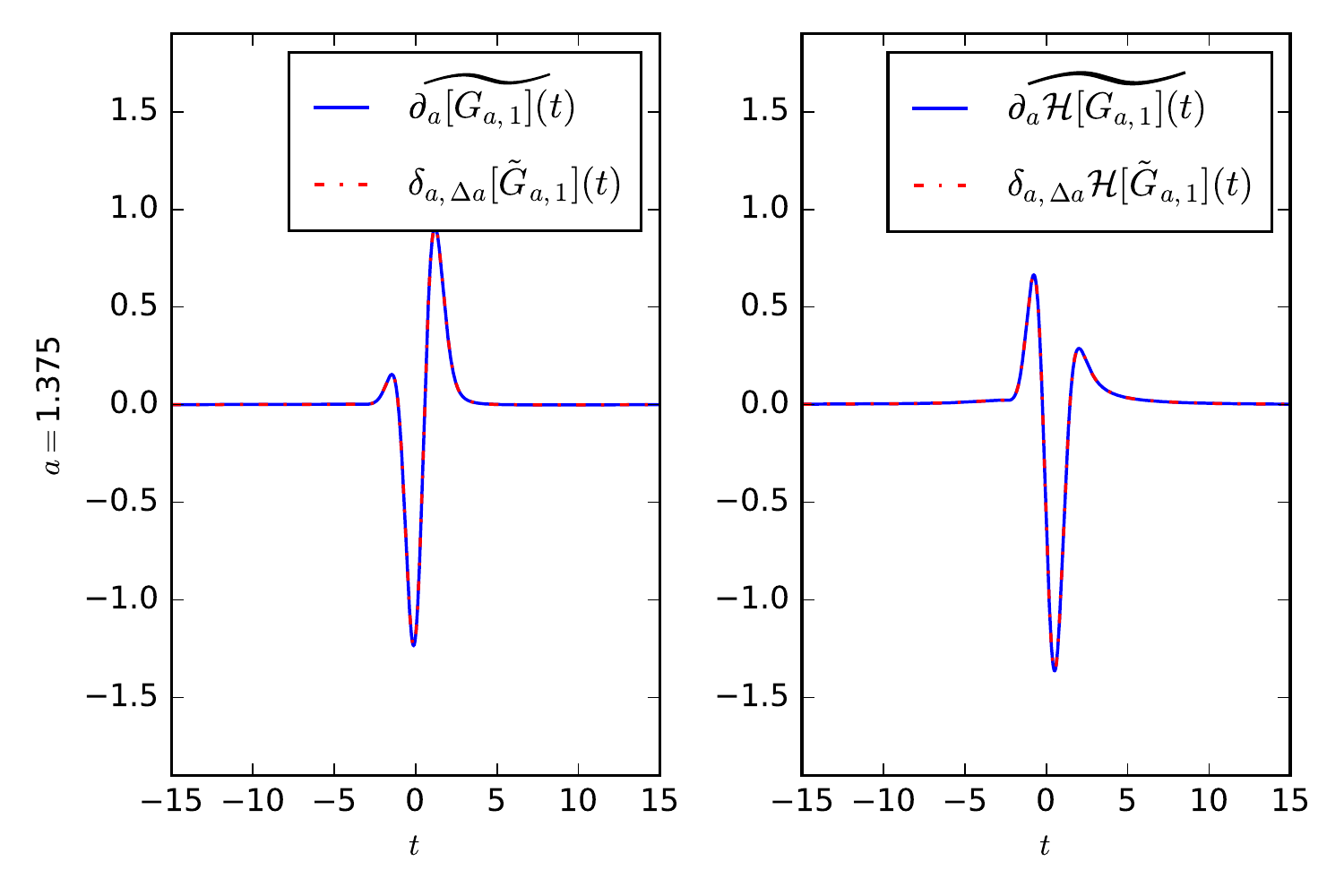} &
\includegraphics[width=0.3\textwidth]{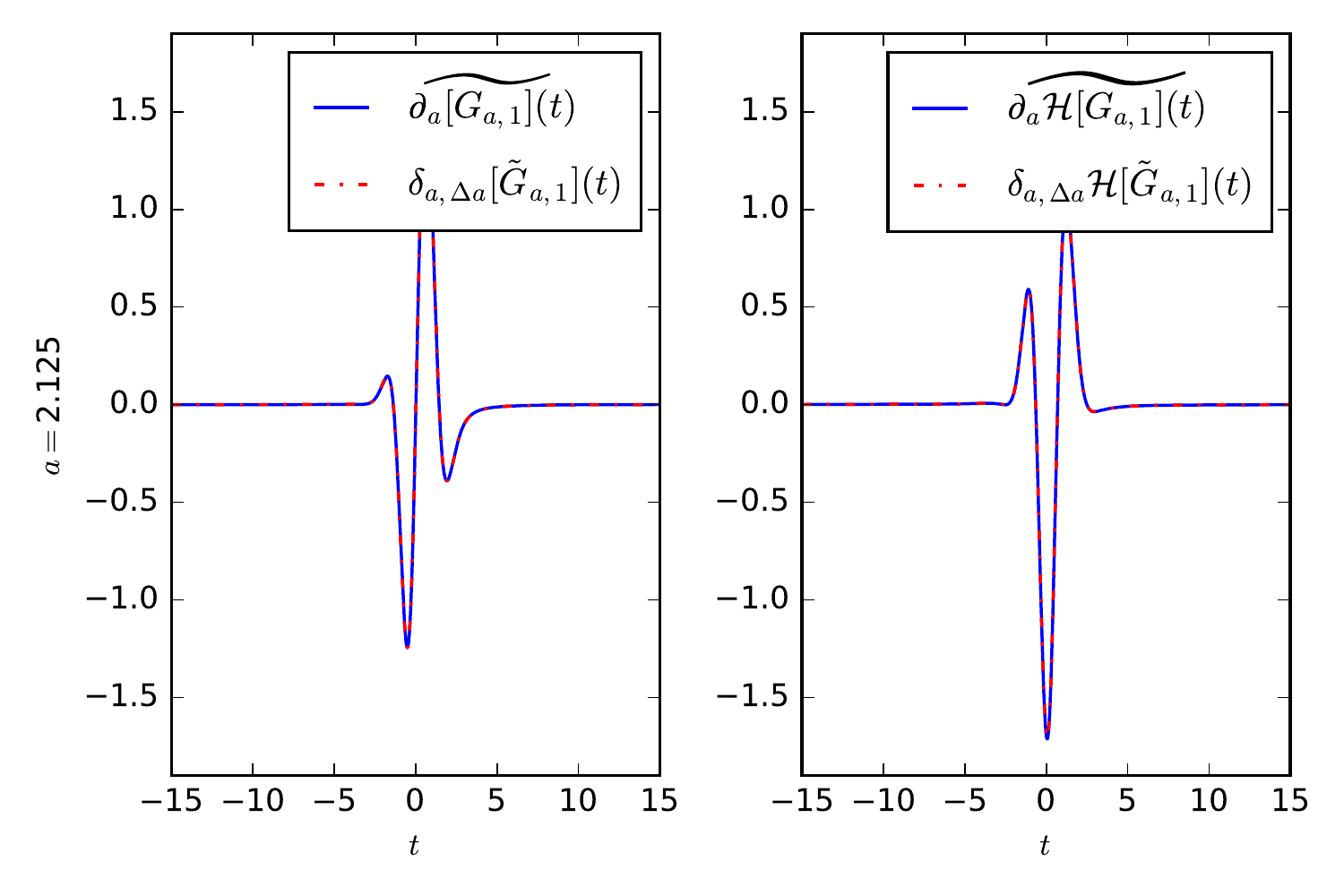}\tabularnewline
\hline 
\end{tabular}\caption{\label{fig:Plot-of-partial_a_of_fractional-derivative} Table of plots
of $\widetilde{\partial_{a}\left[G_{a,\sigma}\right]}\left(t\right)$
(left, solid blue), $\delta_{a,\Delta a}\left[\tilde{G}_{a,\sigma}\right]\left(t\right)$
(left, dashed-dotted red), and $\widetilde{\partial_{a}\mathcal{H}\left[G_{a,\sigma}\right]}\left(t\right)$
(right, solid blue), $\delta_{a,\Delta a}\mathcal{H}\left[\tilde{G}_{a,\sigma}\right]\left(t\right)$
(right,dashed-dotted red), for $\sigma=\sqrt{2}^{-1}$, $\Delta a=10^{-3}$
and various $a$ starting from $0$ to $2.125$ with increments of
$0.125$. Along each column of the table, $a$ increases top to bottom
by $0.125$. Along each row of the table, $a$ increases left to right
by $0.75$. }
\end{figure}

We can generalize the approach presented in Section \ref{subsec:Approximating-fractional-derivative}
to approximate partial derivative with respect to the order of fractional
derivative.{\small{} Consider }(\ref{eq:f_a_sigma_taylor}): {\small{}
\begin{align*}
f_{a,\sigma}\left(t\right) & =\frac{1}{\sqrt{\pi}}\exp\left(\mathrm{i}a\frac{\pi}{2}\right)\sum_{n=0}^{\infty}\frac{\mathrm{i}^{n}}{n!}\left(\frac{\sqrt{2}}{\sigma}\right)^{n+a}\Gamma\left(\frac{a+n+1}{2}\right)t^{n}.
\end{align*}
  Then 
\begin{multline*}
\partial_{a}\left[f_{a,\sigma}\right]\left(t\right)\\
=\frac{1}{\sqrt{\pi}}\exp\left(\mathrm{i}a\frac{\pi}{2}\right)\sum_{n=0}^{\infty}\frac{\mathrm{i}^{n}}{n!}\Gamma\left(\frac{a+n+1}{2}\right)\left(\frac{\sqrt{2}}{\sigma}\right)^{n+a}\left[\log\left(\frac{\sqrt{2}}{\sigma}\right)+\frac{1}{2}\psi\left(\frac{a+n+1}{2}\right)+\mathrm{i}\frac{\pi}{2}\right]t^{n}.
\end{multline*}
Here we used the identity $\Gamma'\left(z\right)=\Gamma\left(z\right)\psi\left(z\right)$
to write the derivative of Gamma function where $\psi\left(z\right)$
is the Psi  function (see 5.2(i), 5.7(ii) and 5.9(ii) of \cite{olver2010nist}):
\begin{align*}
\psi\left(x\right) & =\int_{0}^{\infty}\left(\frac{\mathrm{e}^{-t}}{t}-\frac{\mathrm{e}^{-xt}}{1-\mathrm{e}^{-t}}\right)dt=-\gamma+\sum_{n=0}^{\infty}\left(\frac{1}{n+1}-\frac{1}{n+z}\right).
\end{align*}
Then,} {\small{}partial derivative $\partial_{a}\left[f_{a,\sigma}\right]\left(t\right)$
of $f_{a,\sigma}\left(t\right)$} with respect to the order of fractional
derivative $a$ can be approximated by 
\begin{align*}
\widetilde{\partial_{a}\left[f_{a,\sigma}\right]}\left(t\right) & =\widetilde{\partial_{a}\left[G_{a,\sigma}\right]}\left(t\right)+\mathrm{i}\widetilde{\partial_{a}\mathcal{H}\left[G_{a,\sigma}\right]}\left(t\right)=\sum_{m=1}^{M}\alpha_{m}g\left(\gamma_{m}t\right).
\end{align*}
where $g\left(t\right)$ was defined in (\ref{eq:g(t)}) and $\left(\alpha_{m},\gamma_{m}\right)$
is obtained by solving the following moment problem: 
\begin{multline}
\frac{1}{\sqrt{\pi}}\exp\left(\mathrm{i}a\frac{\pi}{2}\right)\left(\frac{\sqrt{2}}{\sigma}\right)^{a+n}\frac{\Gamma\left(\frac{n+2}{2}\right)\Gamma\left(\frac{a+n+1}{2}\right)\left[\log\left(\frac{\sqrt{2}}{\sigma}\right)+\frac{1}{2}\psi\left(\frac{a+n+1}{2}\right)+\mathrm{i}\frac{\pi}{2}\right]}{n!}\\
=\sum_{m}\alpha_{m}\gamma_{m}^{n}+\varepsilon_{n}.\label{eq:partial_a_frac_derivative-moment_problem}
\end{multline}
 We compared $\widetilde{\partial_{a}\left[f_{a,\sigma}\right]}\left(t\right)$
with the two sided finite difference derivative approximation $\delta_{a,\Delta a}\left[\tilde{f}_{a,\sigma}\right]\left(t\right)$
of $\tilde{f}_{a,\sigma}\left(t\right)=\tilde{G}_{a,\sigma}\left(t\right)+\mathrm{i}\mathcal{H}\left[\tilde{G}_{a,\sigma}\right]\left(t\right)$
with respect to $a$: 
\begin{align*}
\delta_{a,\Delta a}\left[\tilde{f}_{a,\sigma}\right]\left(t\right) & =\frac{\tilde{f}_{a+\Delta a,\sigma}\left(t\right)-\tilde{f}_{a-\Delta a,\sigma}\left(t\right)}{2\Delta a},
\end{align*}
 for some small positive $\Delta a$. For $\Delta a=10^{-3}$, we
present plots of real and imaginary parts of $\widetilde{\partial_{a}\left[f_{a,\sigma}\right]}\left(t\right)$
and $\delta_{a,\Delta a}\left[\tilde{f}_{a,\sigma}\right]\left(t\right)$
for various values of $a$ in Figure \ref{fig:Plot-of-partial_a_of_fractional-derivative}.

\appendix

\section{A rational approximation of Dawson's integral \label{sec:A-rational-approximation2dawson}}

\begin{figure}
\includegraphics[width=0.3\textwidth]{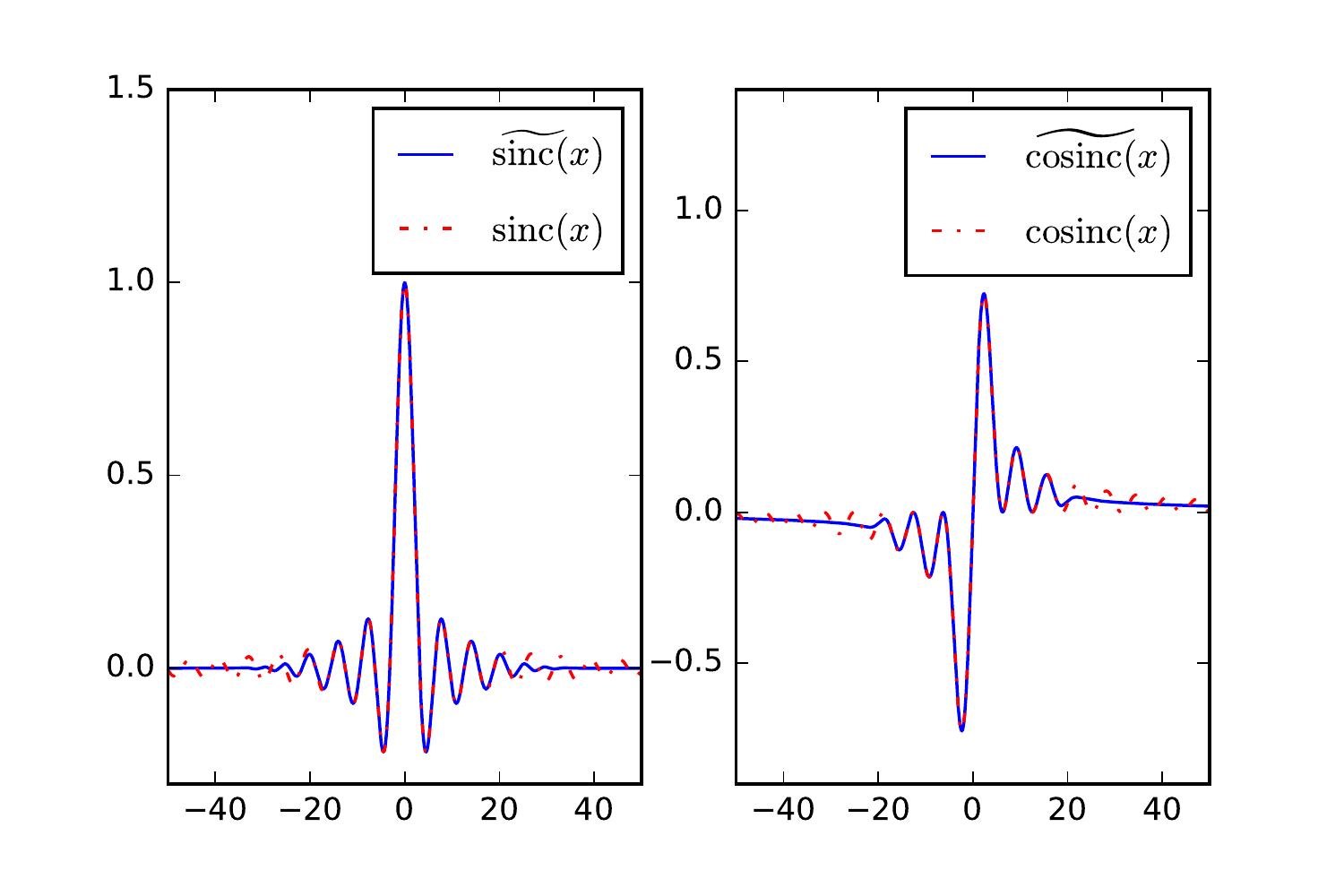}\caption{\label{fig:Approximation-of-sinc-in-Gaussian} Plots of approximations
of $\widetilde{\mbox{sinc}}\left(x\right)$ (left, solid blue) and
$\widetilde{\mbox{cosinc}}\left(x\right)$ (right, solid blue), and
$\mbox{sinc}\left(x\right)$ (left, dashed-dotted  red) and $\mbox{cosinc}\left(x\right)$
(right, dashed-dotted  red).}

\end{figure}

Consider $\text{sinc}\left(x\right)=x^{-1}\sin\left(x\right)$ ,$\mbox{cosinc}\left(x\right)=x^{-1}\left(1-\cos\left(x\right)\right)$
and the approximation 

\begin{align}
\text{sinc}\left(x\right)+\mathrm{i}\text{cosinc}\left(x\right) & =\sum_{m}\alpha_{m}\left[\exp\left(-\left[\gamma_{m}x\right]^{2}\right)+\mathrm{i}\frac{2}{\sqrt{\pi}}F\left(\gamma_{m}x\right)\right]+\epsilon_{B}\left(x\right)\label{eq:epsilon_B-1}
\end{align}
with error $\epsilon_{B}\left(x\right)$, such that $\mathrm{Re}\left\{ \gamma_{m}\right\} >0$.
Here $\left(\alpha_{m},\gamma_{m}\right)$ is a solution to the moment
problem 
\begin{align}
\frac{\Gamma\left(\frac{n+2}{2}\right)}{\left(n+1\right)!} & =\sum_{m}\alpha_{m}\gamma_{m}^{n}+\epsilon_{n},\label{eq:sinc-in-gaussian-moment-problem}
\end{align}
for some $\epsilon_{n}$, which can be solved using the method in
\ref{eq:sinc-in-gaussian-moment-problem}. A solution for the moment
problem \ref{eq:sinc-in-gaussian-moment-problem} is presented in
Table \ref{tab:(alpha_m,gamma_m)_sinc+1i*cosinc_in_gaussian+1i*dawson}
and plots of the corresponding approximations 
\begin{align}
\widetilde{\mbox{sinc}}\left(x\right) & =\mathrm{Re}\left\{ \sum_{m}\alpha_{m}\left[\exp\left(-\left[\gamma_{m}x\right]^{2}\right)+\mathrm{i}\frac{2}{\sqrt{\pi}}F\left(\gamma_{m}x\right)\right]\right\} \nonumber \\
\widetilde{\mbox{cosinc}}\left(x\right) & =\mathrm{Imag}\left\{ \sum_{m}\alpha_{m}\left[\exp\left(-\left[\gamma_{m}x\right]^{2}\right)+\mathrm{i}\frac{2}{\sqrt{\pi}}F\left(\gamma_{m}x\right)\right]\right\} \label{eq:approx-sinc-cosinc}
\end{align}
 are presented in Figure \ref{fig:Approximation-of-sinc-in-Gaussian}.

\begin{table}
\begin{centering}
\begin{tabular}{|c|c|c|}
\hline 
$m$  &
$\alpha_{m}$  &
$\gamma_{m}$\tabularnewline
\hline 
\hline 
{\footnotesize{}1 } &
{\footnotesize{}-0.117532571756027 - 0.003575367485193i} &
{\footnotesize{}0.003599251866768 + 0.023031988828263i}\tabularnewline
\hline 
{\footnotesize{}2 } &
{\footnotesize{}-0.117532571756027 + 0.003575367485193i} &
{\footnotesize{}0.003599251866768 - 0.023031988828263i}\tabularnewline
\hline 
{\footnotesize{}3 } &
{\footnotesize{}2.658250413824904 - 0.872108421067261i} &
{\footnotesize{}0.015512727403264 + 0.023419328488677i}\tabularnewline
\hline 
{\footnotesize{}4 } &
{\footnotesize{}2.658250413824904 + 0.872108421067261i} &
{\footnotesize{}0.015512727403264 - 0.023419328488677i}\tabularnewline
\hline 
{\footnotesize{}5 } &
{\footnotesize{}9.141742362072840 - 34.331838882055756i} &
{\footnotesize{}0.032137996250724 + 0.006297276091330i}\tabularnewline
\hline 
{\footnotesize{}6 } &
{\footnotesize{}9.141742362072840 + 34.331838882055756i} &
{\footnotesize{}0.032137996250724 - 0.006297276091330i}\tabularnewline
\hline 
{\footnotesize{}7 } &
{\footnotesize{}-11.182460204141808 + 11.443034213692357i} &
{\footnotesize{}0.026033590315800 + 0.017203623738202i}\tabularnewline
\hline 
{\footnotesize{}8 } &
{\footnotesize{}-11.182460204141808 - 11.443034213692357i} &
{\footnotesize{}0.026033590315800 - 0.017203623738202i}\tabularnewline
\hline 
\end{tabular}
\par\end{centering}
\caption{\label{tab:(alpha_m,gamma_m)_sinc+1i*cosinc_in_gaussian+1i*dawson}
Table of $(\alpha_{m},\gamma_{m})$ obtained by solving the moment
problem (\ref{eq:sinc-in-gaussian-moment-problem}) using the method
of (\ref{eq:sinc-in-gaussian-moment-problem}). These $(\alpha_{m},\gamma_{m})$
are used in (\ref{eq:approx-sinc-cosinc}) to approximate $\mbox{sinc}\left(x\right)+\mathrm{i}\,\mbox{cosinc}\left(x\right)$
and in (\ref{eq:approx_Dawson}) to approximate Dawson's integral.}
\end{table}

Rewriting Dawson's integral (\ref{eq:dawson-definition}),  
\begin{align*}
F\left(x\right) & =\frac{x}{2}\int_{0}^{\infty}\mathrm{e}^{-t^{2}/4}\text{sinc}\left(xt\right)t\,dt
\end{align*}
 and substituting the approximation (\ref{eq:epsilon_B-1}) for $\text{sinc}\left(x\right)$,
we approximate Dawson's integral by 
\begin{align}
\tilde{F}\left(x\right) & =\frac{x}{2}\sum_{m}\alpha_{m}\int_{0}^{\infty}\mathrm{e}^{-\left(1/4+\left(\gamma_{m}x\right)^{2}\right)t^{2}}t\,dt=x\sum_{m}\frac{\alpha_{m}}{\left(1+\left[2\gamma_{m}x\right]^{2}\right)}.\label{eq:approx_Dawson}
\end{align}

\subsection*{Approximation error}

For $x$ within the vicinity of zero we have the error  
\begin{align*}
\epsilon_{F}\left(x\right) & =\left[F-\tilde{F}\right]\left(x\right)=\frac{x}{2}\int_{0}^{\infty}\mathrm{e}^{-t^{2}/4}\left[\text{sinc}\left(xt\right)-\sum_{m}\alpha_{m}\mathrm{e}^{-\left(\gamma_{m}xt\right)^{2}}\right]t\,dt
\end{align*}
 absolutely bounded by 
\begin{align*}
\left|\epsilon_{F}\left(x\right)\right| & \le\frac{1}{2}\int_{0}^{\infty}\mathrm{e}^{-\left(t/x\right)^{2}/4}\frac{t}{\left|x\right|}\left|\text{sinc}\left(t\right)-\sum_{m}\alpha_{m}\mathrm{e}^{-\left(\gamma_{m}t\right)^{2}}\right|\,dt\\
 & \le\frac{\epsilon_{1}}{2}\int_{0}^{\infty}\mathrm{e}^{-\left(t/x\right)^{2}/4}\frac{t}{\left|x\right|}dt=\frac{\epsilon_{1}}{2}\left|x\right|,
\end{align*}
 where 
\begin{align*}
\epsilon_{1} & =\max_{x\in\mathbb{R}}\left|\text{Re}\left\{ \epsilon_{B}\left(x\right)\right\} \right|,
\end{align*}
 for $B=1$.

For $x$ away from zero  
\begin{align*}
\left|\epsilon_{F}\left(x\right)\right| & \le\frac{1}{2x}\int_{0}^{\infty}\mathrm{e}^{-\left(t/x\right)^{2}/4}\left|\sin t-t\sum_{m}\alpha_{m}\mathrm{e}^{-\left(\gamma_{m}t\right)^{2}}\right|\,dt\\
 & \le\frac{1}{2x}\int_{0}^{\infty}\mathrm{e}^{-\left(t/x\right)^{2}/4}\,dt\left[1+\max_{t\in\mathbb{R}^{+}}\left|t\sum_{m}\alpha_{m}\mathrm{e}^{-\left(\gamma_{m}t\right)^{2}}\right|\right]\\
 & \le\frac{1}{2}\sqrt{\frac{\pi}{x}}\left[1+\sum_{m}\frac{\left|\alpha_{m}\right|}{\sqrt{2\text{Re}\left\{ \gamma_{m}^{2}\right\} }}\mathrm{e}^{-1/2}\right].
\end{align*}
 Here, we used the inequality 
\begin{align*}
\max_{t\in\mathbb{R}^{+}}\left\{ \left|t\sum_{m}\alpha_{m}\mathrm{e}^{-\left(\gamma_{m}t\right)^{2}}\right|\right\}  & \le\sum_{m}\left|\alpha_{m}\right|\max_{t\in\mathbb{R}^{+}}\left\{ t\mathrm{e}^{-\text{Re}\left\{ \gamma_{m}^{2}\right\} t^{2}}\right\} \\
 & =\sum_{m}\frac{\left|\alpha_{m}\right|}{\sqrt{2\text{Re}\left\{ \gamma_{m}^{2}\right\} }}\mathrm{e}^{-1/2}
\end{align*}
 and obtained using $\max_{t\in\mathbb{R}^{+}}t\mathrm{e}^{-at^{2}}=\sqrt{\left(2a\right)^{-1}}\mathrm{e}^{-1/2}$
 .

\newpage{}

\bibliographystyle{amsplain}
\bibliography{37D__Documents_MATLAB_Beams__1+2_D-kernels_Refe___nction_references_on_sine_integral_function,38D__Documents_MATLAB_Beams_FBMR_references_analytic_rep,39D__Schlumberger_LaTeX_Notes_Bandwidth_Extention_slb-bib}

{\small{}}{\small \par}
\end{document}